\documentclass[a4paper, 12pt,reqno]{amsart}

\usepackage[margin=1in]{geometry, stmaryrd}
\usepackage{dirtytalk}
\usepackage{bigints}
\usepackage{amssymb,latexsym}
\usepackage{graphicx}
\usepackage{float}
\graphicspath{ {./images/} }
\usepackage{mathtools}
\usepackage{color}
\usepackage[hyphens]{url}
\usepackage[T1]{fontenc}
\usepackage{hyperref}
\usepackage{amsmath}
\usepackage{mathdots}
\usepackage{breqn}
\usepackage[toc,page]{appendix}
\usepackage{url}    
\usepackage{breqn}
\usepackage{breakurl}
\usepackage{comment}
\usepackage{thmtools}
\usepackage{thm-restate}
\usepackage{mathrsfs}
\usepackage{fancyhdr} 
\usepackage{placeins}
\usepackage{leftidx}
\usepackage{bbm, dsfont}
\allowdisplaybreaks
\makeatletter
\def\imod#1{\allowbreak\mkern10mu({\operator@font mod}\,\,#1)}
\makeatother
\theoremstyle{plain}
\numberwithin{equation}{section}
\newtheorem{thm}{Theorem}[section]
\newtheorem*{theorem*}{Theorem}
\newtheorem{theorem}[thm]{Theorem}
\newtheorem{lemma}[thm]{Lemma}

\newtheorem{definition}[thm]{Definition}
\newtheorem{question}[thm]{Question}
\newtheorem{proposition}[thm]{Proposition}

\newtheorem{remark}[thm]{Remark}

\newcommand{\sfrac}{\genfrac{}{}{}1}
\newcommand{\bigslant}[2]{{\raisebox{.2em}{$#1$}\left/\raisebox{-.2em}{$#2$}\right.}}
\DeclareMathOperator{\GL}{GL}
\DeclareMathOperator{\SO}{SO}
\DeclareMathOperator{\Su}{Supp}
\DeclareMathOperator{\Hom}{Hom}
\DeclareMathOperator{\N}{N}
\DeclareMathOperator{\F}{F}
\DeclareMathOperator{\E}{E}

\DeclareMathOperator{\W}{W}
\DeclareMathOperator{\Mat}{Mat}
\DeclareMathOperator{\Ind}{Ind}
\DeclareMathOperator{\ind}{ind}
\DeclareMathOperator{\SP}{span}

\title{O\lowercase{n} \lowercase{the} L\lowercase{ocal} L\lowercase{anglands} F\lowercase{unctorial} T\lowercase{ransfer} \lowercase{from} $\SO(5)$ \lowercase{to} $\GL(4)$}
\author{D\lowercase{avid} C.~L\lowercase{uo}}
\date{}

\address{School of Mathematics, University of Minnesota, Minneapolis, MN 55455, United States} \email{luo00275@umn.edu}

\subjclass{22E50, 11F70}
\keywords{local Langlands functorial transfer, local Shalika model, supercuspidal representation, local exterior square $L$-function.}

\begin{document}

\begin{abstract}
    Let $\F$ be a non-archimedean local field of characteristic zero.~We study the local Langlands functorial transfer from $\SO(5, \F)$ to $\GL(4, \F)$ for several families of supercuspidal representations of $\GL(4, \F)$.~In doing so, we give equivalent conditions for such representations to be functorial transfers from $\SO(5, \F)$, expressed in terms of the Bushnell--Kutzko construction of supercuspidal representations, by studying poles of local exterior square $L$-functions and the existence of non-zero local Shalika models.~
\end{abstract}

\maketitle


\section{Introduction}

Let $\F$ be a non-archimedean local field of characteristic zero and $\pi$ a supercuspidal representation of the degree four general linear group $\GL(4, \F)$.~An important invariant one may attach to special $\pi$ is a non-zero \textit{local Shalika model}.~The existence of non-zero local Shalika models enables us to detect local Langlands functorial transfers from the degree five special orthogonal group $\SO(5, \F)$ and provides various other applications to the theory of automorphic forms in the global setting \cite[Theorem 5.5, Section 6]{Jiang}.~In the local setting, the principle of functoriality gives a conceptual framework for transferring admissible representations of general split reductive algebraic groups over $\F$ \cite[Lecture 11]{Cogdell}.

In 2008, Jiang, Nien, and Qin proved that a supercuspidal representation $\tau$ of $\GL(2n, \F)$ is a functorial transfer from $\SO(2n+1, \F)$ if and only if $\tau$ possesses a non-zero local Shalika model \cite[Theorem 5.5]{Jiang}.~To show that $\tau$ possesses such a model, Jiang, Nien, and Qin further equated this condition to the local exterior square $L$-function $L(s, \tau, \wedge^{2})$ having a pole at $s = 0$ \cite[Theorem 5.5]{Jiang}.~For supercuspidals of $p$-adic general linear groups, $L(s, \tau, \wedge^{2})$ is well understood and can be written explicitly as a product in terms of the function $\Lambda_{s_{0}}$ that we define in Subsection \ref{TwistedShalikaPeriod} \cite[Theorem 6.2 (i)]{Jo1}, \cite[Definition 3.18]{YZ1}.

The aim of this article is to shed light on the following question:
\begin{question}\label{Question}
    Given a supercuspidal representation of $\GL(4, \F)$ that possesses a non-zero local Shalika model, what can be deduced about its construction in terms of type theory?
\end{question}
\noindent Answering Question \ref{Question} in the setting of $\GL(4, \F)$ would provide a starting point for obtaining a more explicit description of the local Langlands functorial transfer from $\SO(2n+1, \F)$ to $\GL(2n, \F)$ for arbitrary $n$.~

In this paper, we consider depth-zero and three classes of minimax (in the sense of \cite{Adrian}) supercuspidal representations of $\GL(4, \F)$.~The minimax supercuspidals we consider are \textit{simple} (\textit{epipelagic}) \textit{supercuspidal representations} (those minimax supercuspidals of minimal positive depth $\frac{1}{4}$) \cite{Knightly, Adrian-Liu}, \textit{middle supercuspidal representations} (those minimax supercuspidals of next smallest positive depth $\frac{1}{2}$) \cite{Luo}, and \textit{biquadratic supercuspidal representations} (those minimax supercuspidals of depth-one which we introduce in Subsection \ref{Biquadratic}).~Such minimax supercuspidal representations are induced from simple strata of the form $\left[\mathfrak{A}, 1, 0, \beta\right]$ whose associated set of simple characters is a singleton in the Bushnell--Kutzko sense \cite{BK}.

In the cases of depth-zero, middle, and biquadratic supercuspidal representations, the requirement for such supercuspidals to possess non-zero local Shalika models is minimal.~For a simple supercuspidal representation $\pi$, this requirement becomes more complicated.~To understand the necessary conditions for $\pi$ to possess a non-zero local Shalika model, we first parametrize $\pi$ by the triple $(\overline{v}, \, \varphi, \, \zeta)$ where $\overline{v}$ is a non-zero element of the residue field $k_{\F}$ of $\F$, $\varphi$ is a quasi-character of $k_{\F}^{\times}$, and $\zeta \in \mathbb{C}^{\times}$.~We may use such a parametrization as there exists a bijection between the set of isomorphism classes of simple supercuspidal representations and the set of all such triples \cite[Proposition 1.3]{Imai}.~

Using the above classification, we obtain the following main result.~The cases of depth-zero and simple supercuspidal representations are due to Ye and Zelingher in the respective references: \cite[Theorems 3.14 and 3.15]{YZ1}, \cite[Theorem 3.5 (2)]{YZ2}.
\begin{theorem}\label{Main}
    Let $\pi$ be a depth-zero, middle, or biquadratic supercuspidal representation of $\GL(4, \F)$.~The following are equivalent:
    \begin{enumerate}
        \item $\pi$ is a local Langlands functorial transfer from $\SO(5, \F)$.
        \item $\pi$ possesses a non-zero local Shalika model.~
        \item $\pi$ has trivial central character.
    \end{enumerate}
    Let $\pi := \pi_{(\overline{v}, \, \varphi, \, \zeta)}$ be a simple supercuspidal representation of $\GL(4, \F)$.~The following are equivalent:
    \begin{enumerate}
        \item $\pi$ is a local Langlands functorial transfer from $\SO(5, \F)$.
        \item $\pi$ possesses a non-zero local Shalika model.~
        \item $\pi$ has trivial central character and $\zeta = \pm 1$.
    \end{enumerate}
\end{theorem}
If we assume that $\F$ has odd residual characteristic, then Theorem \ref{Main} also holds for supercuspidal representations of the form $\pi \otimes \Psi$, where $\pi$ is a depth-zero, simple, middle, or biquadratic supercuspidal representation, and $\Psi$ is a quasi-character of $\GL(4, \F)$ (see Remark \ref{RemarkMain} for more details).~

For positive depth minimax supercuspidals, we seek to understand what separates the cases of middle and biquadratic supercuspidals from the case of simple supercuspidals in Theorem \ref{Main}.~To begin, we distinguish these representations via their \textit{minimal polynomials} (see Definition \ref{Characteristic}).~For a general supercuspidal representation of $\GL(4, \F)$, its minimal polynomial is an invariant attached to its defining simple stratum.~From \cite[Lemma 4.2]{Luo}, we see that the minimal polynomial attached to a simple supercuspidal is linear, whereas for a middle supercuspidal, it is quadratic.~Furthermore, we have from Subsection \ref{Biquadratic} that the minimal polynomial attached to a biquadratic supercuspidal is biquadratic, i.e.~a polynomial omitting all odd power terms.~

Outside of biquadratic supercuspidals, one may also construct other depth-one minimax supercuspidals induced from the maximal parahoric subgroup whose associated set of simple characters is a singleton.~Namely, these are the supercuspidals whose minimal polynomials include odd power terms.~We note, however, that the techniques we use for middle and biquadratic supercuspidals do not apply to these alternative depth-one minimax supercuspidals because $\Lambda_{0}$, when evaluated at our chosen Whittaker function, is zero.~In principle, there could be a Whittaker function which gives a non-zero $\Lambda_{0}$ value.~On the other hand, it would be interesting to know the following:
\begin{question}
    Let $\pi$ be an irreducible minimax supercuspidal representation of $\GL(4, \F)$. If the minimal polynomial associated with $\pi$ is quadratic or biquadratic, is it true that $\pi$ possesses a non-zero local Shalika model if and only if $\pi$ has trivial central character?
\end{question}

The structure of this paper is given as follows.~In Section \ref{PreliminariesAndBackground}, we provide the necessary background information needed for this article.~Specifically, we give a brief account of maximal simple types and the Bushnell--Kutzko construction of supercuspidal representations of $\GL(4, \F)$.~Next, we introduce an explicit Whittaker function given by Pa\v{s}k\={u}nas and Stevens which will enable us to carry out our computations involving $\Lambda_{s_{0}}$.

In Subsection \ref{FuncTrans}, we define what it means for an irreducible admissible representation of $\GL(4, \F)$ to be a local Langlands functorial transfer from $\SO(5, \F)$.~Then in Subsection \ref{TwistedShalikaPeriod}, we discuss how the local exterior square $L$-function for a supercuspidal representation of $\GL(4, \F)$ is related to $\Lambda_{s_{0}}$ and give a general strategy for proving Theorem \ref{Main}.~

In Section \ref{ConstructionSupercuspidals}, we first recall the construction of simple and middle supercuspidal representations and explicit Pa\v{s}k\={u}nas--Stevens Whittaker functions in their Whittaker models.~Next, we construct biquadratic supercuspidal representations and give the associated Pa\v{s}k\={u}nas--Stevens Whittaker function.~In Section \ref{ProofOfTheorem}, we first show that middle and biquadratic supercuspidal representations possess non-zero local Shalika models and are local Langlands functorial transfers from $\SO(5, \F)$ if and only if their central characters are trivial.~Finally, we give a proof of Theorem \ref{Main} in Subsection \ref{End}.

\section*{Acknowledgments}
The author is grateful to Dihua Jiang for suggesting this problem and to Shaun Stevens and Rongqing Ye for valuable discussions.

\section{Preliminaries and background}\label{PreliminariesAndBackground}

In this section, we give the necessary preliminary and background information needed for the remainder of this article.

\subsection{Notation}\label{Notation}
To fix some notation, let $\F$ be a non-archimedean local field of characteristic zero with $\mathcal{O}_{\F}$ its valuation ring, $k_{\F}$ its residue field with cardinality $q_{\F}$ and characteristic $p$, and $\mu_{\F}'$ the group of roots of unity of order prime to $p$.~Next, let $\mathcal{P}_{\F}$ denote the maximal ideal of $\mathcal{O}_{\F}$ and $\varpi_{\F}$ a fixed uniformizer.~Furthermore, let $\psi_{\F}$ be a non-trivial additive quasi-character of $\F$ with conductor $\mathcal{P}_{\F}$, i.e.~$\psi_{\F}$ is non-trivial on $\mathcal{O}_{\F}$ but trivial on $\mathcal{P}_{\F}$.~ 

Let $\N(n, \F)$ denote the unipotent radical of the standard Borel subgroup of invertible upper triangular matrices of $\GL(n, \F)$, $\text{P}(n, \F)$ the standard mirabolic subgroup, $\mathcal{B}(n, \F)$ the subspace of upper triangular $n \times n$ matrices, $\Mat(n \times m, R)$ the space of $n \times m$ matrices with coefficients in some ring $R$, and $I_{n}$ the $n \times n$ identity matrix of $\GL(n, \F)$.~Next, we let $\psi_{n}$ denote the standard smooth non-degenerate quasi-character of $\N(n, \F)$ defined by
\[
\psi_{n}(u) = \psi_{\F}\left(\sum_{i=1}^{n-1}u_{i, i+1}\right)
\]
where $u = (u_{i, j}) \in \N(n, \F)$.~Moreover, let $\sigma_{4}$ denote the permutation matrix $(23) \in S_{4}$.~

Let $\Ind$ denote smooth induction and $\ind$ compact induction.~Lastly, we will fix the standard normalizations for the Haar measures $dx$ and $d^{\times}x$ on $\F$ and $\F^{\times}$ relative to $\psi_{\F}$ respectively.~In particular, this means that $\int_{\mathcal{O}_{\F}}  dx = 1$ and $\int_{\mathcal{O}_{\F}^{\times}}  d^{\times}x = 1$.

\subsection{Maximal simple types}\label{MaximalSimpleTypes}

From \cite{BK}, we have that supercuspidal representations of $\GL(4, \F)$ may be classified in terms of \textit{maximal simple types} which are special pairs $(J,  \lambda)$ where $J$ is a compact open subgroup of $\GL(4, \F)$
and $\lambda$ is an irreducible representation of $J$.~We refer to \cite{BK} for precise definitions of, and results concerning, the objects introduced in this subsection.

Let $V = \F^{4}$ be a four-dimensional vector space over $\F$ with standard basis.~Thus we identify $\text{Aut}_{\F}(V)$ with $\GL(4, \F)$  and $A = \text{End}_{\F}(V)$ with $\Mat(4 \times 4, \F)$.~Let $\mathfrak{A}$ be a principal hereditary $\mathcal{O}_{\F}$-order in $A$ with period $e\left(\mathfrak{A}\right)$ and Jacobson radical $\mathfrak{P}$.~Define
\[
U^{0}\left(\mathfrak{A}\right) = U\left(\mathfrak{A}\right) = \mathfrak{A}^{\times}, \; U^{n}\left(\mathfrak{A}\right) = I_{4} + \mathfrak{P}^{n}, \; n \geq 1.
\]
For $n \geq 0$, choose $\beta \in A$ such that
\[
\beta \in \mathfrak{P}^{-n} \setminus \mathfrak{P}^{1-n},
\]
where $\E = \F \left[\beta\right]$ is a field extension of $\F$ and $\E^{\times}$ normalizes $\mathfrak{A}$.~Provided an additional technical condition is satisfied (namely that the critical exponent $k_{0}\left(\beta, \mathfrak{A}\right) < 0$), these data give a principal simple stratum $\left[\mathfrak{A}, n, 0, \beta \right]$ of $A$.~Take $J := J\left(\beta,\mathfrak{A}\right)$, $J^{1} := J^{1}\left(\beta,\mathfrak{A}\right)$, and $H^{1} := H^{1}\left(\beta,\mathfrak{A}\right)$ where
\[
H^{1} \subset J^{1} \subset J.
\] 
Denote by $\mathcal{C}\left(\mathfrak{A}, \beta, \psi_{\F}\right)$ the set of simple (linear) characters of $H^{1}$.

Recall the following definition of maximal simple types.
\begin{definition}
    \normalfont
The pair $(J,  \lambda)$, where $J$ is a compact open subgroup of $\GL(4, \F)$ and $\lambda$ is an irreducible representation of $J$, is called a \textit{maximal simple type} if one of the following holds:
\begin{enumerate}
    \item $J = J\left(\beta,\mathfrak{A}\right)$ is a compact open subgroup associated to a simple stratum $\left[\mathfrak{A}, n, 0, \beta \right]$ of $A$ as above, such that, if $\E = \F \left[\beta\right]$ and $B = \text{End}_{\E}(V)$, then $\mathfrak{B} = \mathfrak{A} \cap B$  is a maximal $\mathcal{O}_{\E}$-order in $B$.~Moreover, there exists a simple character $\theta \in \mathcal{C}\left(\mathfrak{A}, \beta, \psi_{\F}\right)$ such that
    \[
    \lambda \cong \kappa \otimes \sigma
    \]
    where $\kappa$ is a $\beta$-extension of the unique irreducible representation $\eta$ of $J^{1} = J^{1}\left(\beta,\mathfrak{A}\right)$, which contains $\theta$, and $\sigma$ is the inflation to $J$ of an irreducible cuspidal representation of
    \[
    \bigslant{J}{J^{1}} \cong \bigslant{U\left( \mathfrak{B} \right)}{U^{1}\left( \mathfrak{B} \right)} \cong \GL(r, k_{\E}),
    \]
    where $r =\frac{4}{[\E \, : \, \F]}$.
    
    \item $(J,  \lambda) = \left(U\left( \mathfrak{A} \right),  \sigma \right)$, where $\mathfrak{A}$ is a maximal hereditary $\mathcal{O}_{\F}$-order in $A$ and $\sigma$ is the inflation to $U\left( \mathfrak{A} \right)$ of an irreducible cuspidal representation of
    \[
    \bigslant{U\left( \mathfrak{A} \right)}{U^{1}\left( \mathfrak{A} \right)} \cong \GL(4, k_{\F}).
    \]
\end{enumerate}
Since the condition $\beta \in \mathfrak{P}^{-n} \setminus \mathfrak{P}^{1-n}$ in case (1) implies that $\beta$ is non-zero, the depth-zero situation in case (2) is not included in the preceding construction.~For notational convenience, we will nevertheless regard case (2) formally as a special case of case (1) by setting $\beta = 0$, $\E = \F$, and $\theta$, $\eta$, $\kappa$ trivial.~In either case, we write $\textbf{J} = \E^{\times}J$.~With these data, any supercuspidal representation $\tau$ of $\GL(4, \F)$ is of the form
\[
\tau \cong \ind_{\textbf{J}}^{\GL(4, \, \F)}\Lambda
\]
for some choice of $(\textbf{J},  \Lambda)$, where $\Lambda \big|_{J} = \lambda$.~We call such a pair $(\textbf{J},  \Lambda)$ an \textit{extended maximal simple type}.~Let $d(\tau) = \frac{n}{e\left(\mathfrak{A}\right)}$ denote the \textit{depth} of $\tau$.~Lastly, a supercuspidal representation induced from a maximal simple type of the form in case (2) is called a \textit{depth-zero supercuspidal representation}.
\end{definition}
\noindent For $\tau$ a supercuspidal representation of $\GL(4, \F)$, any two extended maximal simple types in $\tau$ are conjugate in $\GL(4, \F)$.~

Finally, we recall that for a positive depth supercuspidal representation $\tau$ of $\GL(4, \F)$, its associated simple stratum $\left[\mathfrak{A}, n, 0, \beta\right]$ has a \textit{characteristic polynomial} $\varphi_{Y_{\beta}}(X) \in k_{\F}[X]$ that depends only on the equivalence class of the final simple stratum in the defining sequence of $\left[\mathfrak{A}, n, 0, \beta\right]$ \cite[Chapter 2, p.~58]{BK}.~
\begin{definition}\label{Characteristic}
    \normalfont
    Let $e := e\left(\mathfrak{A}\right)$ and $g := \gcd(e, n)$.~We set 
    \[
    Y_{\beta} = \beta^{e/g}\varpi_{\F}^{n/g} \in A.
    \]
    Let $\varphi_{Y_{\beta}}(X) \in \F[X]$ be the characteristic polynomial of $Y_{\beta}$ as an $\F$-endomorphism of $V$.~Then, since $Y_{\beta} \in \mathfrak{A}$, we have $\varphi_{Y_{\beta}}(X) \in \mathcal{O}_{\F}[X]$.~We call the reduction $\overline{\varphi}_{Y_{\beta}}(X)$ of $\varphi_{Y_{\beta}}(X)$ modulo $\mathcal{P}_{\F}$ the \textit{characteristic polynomial} of $\left[\mathfrak{A}, n, 0, \beta\right]$.~Furthermore, $\overline{\varphi}_{Y_{\beta}}(X)$ is a power of a monic irreducible polynomial $\varphi_{\tau}(X) \in k_{\F}[X]$ that we call the \textit{minimal polynomial} of $\left[\mathfrak{A}, n, 0, \beta\right]$.~From \cite[Subsection 6.1]{BHK}, we have that $\varphi_{\tau}$ is an invariant of $\tau$ with $\varphi_{\tau}(X) \neq X$.~
\end{definition}

\subsection{Explicit Whittaker functions}\label{ExplicitWhittakerfunctions}

In this subsection, we introduce an explicit Whittaker function $W_{\tau}$ in the \textit{Whittaker model} $W\left(\tau, \psi_{\F}\right)$ associated with a supercuspidal representation $\tau$ of $\GL(4, \F)$ \cite{Paskunas}.~To do so, we first introduce Bessel functions of supercuspidal representations of $\GL(4, \F)$ constructed by Pa\v{s}k\={u}nas and Stevens.~We recall the basics of these Bessel functions, which rely on the construction theory of supercuspidal representations of $\GL(4, \F)$ in terms of maximal simple types of Bushnell and Kutzko introduced in Subsection \ref{MaximalSimpleTypes}.~We will use the notation from \cite{BK, Paskunas}.

We recall the following (\cite[Section 2]{Adrian-Liu} and \cite{BH-2}):~an irreducible admissible representation $\left(\tau, V_{\tau}\right)$ of $\GL(4, \F)$ is called \textit{generic} if 
\[
\Hom_{\GL(4, \, \F)}\left(\tau,  \Ind_{\N(4, \, \F)}^{\GL(4, \, \F)}\psi_{4}\right) \neq 0.
\]
By the uniqueness of local Whittaker models, this $\Hom$-space is at most one-dimensional. From Frobenius reciprocity,
\[
\Hom_{\GL(4, \, \F)}\left(\tau,  \Ind_{\N(4, \, \F)}^{\GL(4, \, \F)}\psi_{4}\right) \cong \Hom_{\N(4, \, \F)}\left(\tau\big|_{\N(4, \, \F)},  \psi_{4}\right).
\]
Therefore, $\Hom_{\N(4, \, \F)}\left(\tau \big|_{\N(4, \, \F)},  \psi_{4}\right)$ is also at most one-dimensional.~Assume that $\left(\tau, V_{\tau}\right)$ is generic.~Fix a non-zero functional $l \in \Hom_{\N(4, \, \F)}\left(\tau \big|_{\N(4, \, \F)},  \psi_{4}\right)$, which is unique up to scalar.~The Whittaker function attached to a vector $v \in V_{\tau}$ is defined by
\[
W_{v}(g) := l \left(\tau(g)v\right), \; \text{ for all } g \in \GL(4, \F),
\]
so that $W_{v} \in \Ind_{\N(4, \, \F)}^{\GL(4, \, \F)}\psi_{4}$.~The space
\[
W(\tau, \psi_{\F}) := \{ W_{v} : v \in V_{\tau}\}
\]
is called the \textit{Whittaker model} of $\tau$ and $\GL(4, \F)$ acts on it by right translation.~It is easy to see that the Whittaker model of $\tau$ is independent of the choice of the non-zero functional $l$.

Next, we recall from \cite[Section 5]{Paskunas} the general formulation of Bessel functions.~Let $\mathcal{K}$ be an open compact-modulo-centre subgroup of  $\GL(4, \F)$ and let $\mathcal{U} \subset \mathcal{M} \subset \mathcal{K}$ be compact open subgroups of $\mathcal{K}$.~Let $\tau$ be an irreducible smooth representation of $\mathcal{K}$ and let $\Psi$ be a linear character of $\mathcal{U}$.~Take an open normal subgroup $\mathcal{N}$ of $\mathcal{K}$, which is contained in $\text{Ker}(\Psi) \, \cap \, \text{Ker}(\tau)$.~Let $\chi_{\tau}$ be the (trace) character of $\tau$.~The associated Bessel function $\mathcal{J} \, \colon \, \mathcal{K}  \longrightarrow \mathbb{C}$ of $\tau$ is defined by
\[
\mathcal{J}(g) := \left[ \, \mathcal{U} : \mathcal{N} \, \right]^{-1}\displaystyle\sum_{h \in \mathcal{U}/\mathcal{N}}\Psi(h)^{-1}\chi_{\tau}(gh).
\]
This is independent of the choice of $\mathcal{N}$.~The basic properties of this Bessel function which we will need are given below.
\begin{proposition}\cite[Proposition 5.3]{Paskunas}\label{Paskunas5.3}
Assume that the data introduced above satisfy the following:
\begin{itemize}
    \item $\tau \big|_{\mathcal{M}}$ is an irreducible representation of $\mathcal{M}$;
    \item $\tau \big|_{\mathcal{M}} \cong \Ind_{\mathcal{U}}^{\mathcal{M}}\left(\Psi\right)$.
\end{itemize}
Then the Bessel function $\mathcal{J}$ of $\tau$ enjoys the following properties:
\begin{enumerate}
    \item $\mathcal{J}(1)=1$;
    \item $\mathcal{J}(hg) = \mathcal{J}(gh) = \Psi(h)\mathcal{J}(g)$ for all $h \in \mathcal{U}$ and $g \in \mathcal{K}$;
    \item if $\mathcal{J}(g) \neq 0$, then $g$ intertwines $\Psi$; in particular, if $m \in \mathcal{M}$, then $\mathcal{J}(m) \neq 0$ if and only if $m \in \mathcal{U}$.
\end{enumerate}
\end{proposition}

By \cite[Proposition 1.6]{BH-2}, there is an extended maximal simple type $(\textbf{J}, \Lambda)$ in $\tau$ such that
\[
\Hom_{\N(4, \, \F)  \, \cap \,  \textbf{J}}\left(\psi_{4}, \Lambda\right) \neq 0.
\]
Since $\Lambda$ restricts to a multiple of some simple character $\theta \in \mathcal{C}\left(\mathfrak{A}, \beta, \psi_{\F}\right)$, one finds that $\theta(u) = \psi_{4}(u)$ for all $u \in \N(4, \F) \cap H^{1}$.~As in \cite[Definition 4.2]{Paskunas}, one defines a character 
\[
\Psi_{4} \, \colon \, (J \cap \N(4, \F)) \, H^{1} \longrightarrow \mathbb{C}^{\times}
\]
by
\[
\Psi_{4}(uh) := \psi_{4}(u)\theta(h)
\]
for all $u \in J \cap \N(4, \F)$ and $h \in H^{1}$.~By \cite[Theorem 4.4]{Paskunas}, the data
\[
\mathcal{K} = \textbf{J}, \; \tau = \Lambda, \; \mathcal{M} = (J \cap \text{P}(4, \F))J^{1}, \; \mathcal{U} = (J \cap \N(4, \F))H^{1}, \text{ and } \Psi = \Psi_{4}
\]
satisfy the conditions in Proposition \ref{Paskunas5.3} and hence define a Bessel function $\mathcal{J}$.

Now we define a function $W_{\tau} \, \colon \, \GL(4, \F) \longrightarrow \mathbb{C}$ by
\begin{equation}\label{Whittaker}
    W_{\tau}(g):= \begin{cases}
        \psi_{4}(u)\mathcal{J}(j) \; &\text{ if } g=uj \text{ with } u \in \N(4, \F), \; j \in \textbf{J}, \\
        0 \; &\text{ otherwise,}
    \end{cases}
\end{equation}
which is well-defined by Proposition \ref{Paskunas5.3} (2).~Then, by \cite[Theorem 5.8]{Paskunas}, $W_{\tau} \in W(\tau, \psi_{\F})$ is a Whittaker function for $\tau$.~

\subsection{Local Shalika models}\label{ShalikaPeriods}

In this subsection, we introduce the \textit{local Shalika model} of a supercuspidal representation $(\tau, V_{\tau})$ of $\GL(4, \F)$ \cite[Section 2]{Jiang}.~Let $P_{2, 2}(\F) = M_{2, 2} N_{2, 2}$ be the maximal parabolic subgroup of $\GL(4, \F)$, with 
\[
M_{2, 2} = \GL(2, \F) \times \GL(2, \F),
\]
and 
\[
N_{2, 2} = \left\{n(X) = \begin{pmatrix}
    I_{2} & X \\
    & I_{2}
\end{pmatrix}\right\}.
\]
We define a character
\[
\psi_{N_{2, 2}}\left(n(X)\right) = \psi_{\F}\left(\text{tr}(X)\right).
\]
The stabilizer of $\psi_{N_{2, 2}}$ in $M_{2, 2}$ is $\GL(2, \F)^{\Delta}$, the diagonal embedding of $\GL(2, \F)$ into $M_{2, 2}$.

Next, we denote by 
\[
\mathcal{S}(2, \F) = \GL(2, \F)^{\Delta} \rtimes N_{2, 2}
\]
the \textit{Shalika subgroup} of $\GL(4, \F)$.~Denote by $\psi_{\mathcal{S}(2, \, \F)}$ the extension of $\psi_{N_{2, 2}}$ from $N_{2, 2}$ to the Shalika subgroup $\mathcal{S}(2, \F)$, such that $\psi_{\mathcal{S}(2, \, \F)}$ is trivial on $\GL(2, \F)^{\Delta}$.~The \textit{Shalika functional} of $(\tau, V_{\tau})$ is a non-zero functional in the following space
\[
\Hom_{\mathcal{S}(2, \, \F)}\left(V_{\tau}, \psi_{\mathcal{S}(2, \, \F)}\right).
\]
Equivalently, a Shalika functional is a non-trivial functional $f$ on $V_{\tau}$ satisfying
\[
f(\tau(s)v) = \psi_{\mathcal{S}(2, \, \F)}(s)f(v)
\]
for all $s \in \mathcal{S}(2, \F)$ and $v \in V_{\tau}$.~Therefore, $V_{\tau}$ admits a non-trivial embedding into the full induction $\Ind_{\mathcal{S}(2, \, \F)}^{\GL(4, \, \F)}\psi_{\mathcal{S}(2, \, \F)}$, since by Frobenius reciprocity
\[
\Hom_{\mathcal{S}(2, \, \F)}\left(V_{\tau}, \psi_{\mathcal{S}(2, \, \F)}\right) \cong \Hom_{\GL(4, \, \F)}\left(V_{\tau}, \Ind_{\mathcal{S}(2, \, \F)}^{\GL(4, \, \F)}\psi_{\mathcal{S}(2, \, \F)}\right).
\]

By the local uniqueness of the local Shalika model \cite{Jacquet, Nien}, the dimension of the space $\Hom_{\mathcal{S}(2, \, \F)}\left(V_{\tau}, \psi_{\mathcal{S}(2, \, \F)}\right)$ is at most one.~If it is non-zero, we say that $\tau$ has a local Shalika model.~More precisely, if $\mathcal{L}_{\psi_{\mathcal{S}(2, \, \F)}}$ is a non-zero Shalika functional of $(\tau, V_{\tau})$, the \textit{local Shalika model} of $\tau$ consists of all functions of the form
\begin{equation}
    S_{\psi_{\mathcal{S}(2, \, \F)}, v}(g) := \mathcal{L}_{\psi_{\mathcal{S}(2, \, \F)}}(\tau(g)v)
\end{equation}
for all $v \in V_{\tau}$.~It is clear that $S_{\psi_{\mathcal{S}(2, \, \F)}, v}$ belongs to the space $\Ind_{\mathcal{S}(2, \, \F)}^{\GL(4, \, \F)}\psi_{\mathcal{S}(2, \, \F)}$.

\subsection{Local Langlands functorial transfers}\label{FuncTrans}

In this subsection, we define what it means for a supercuspidal representation $\tau$ of $\GL(4, \F)$ to be a local Langlands functorial transfer from $\SO(5, \F)$.~To begin, we recall the local Langlands correspondence for $\GL(4, \F)$ \cite{Henniart1, Harris}, \cite[Section 5]{Jiang}.~Let $\text{SL}(2, \mathbb{C})$ denote the special linear group of degree two with coefficients in $\mathbb{C}$.
\begin{definition}
    \normalfont
    Let $\W_{\F}$ be the Weil group of $\F$.~A \textit{local Langlands parameter} for $\GL(4, \F)$ is a group homomorphism
    \[
    \varphi \, \colon \, \W_{\F} \times \text{SL}(2, \mathbb{C}) \longrightarrow \GL(4, \mathbb{C})
    \]
    such that the restriction of $\varphi$ to $\W_{\F}$ is continuous with respect to the topology of $\W_{\F}$ and the topology of the complex Lie group $\GL(4, \mathbb{C})$ and the restriction of $\varphi$ to $\text{SL}(2, \mathbb{C})$ is algebraic.~
\end{definition}

From the local Langlands correspondence for $\GL(4, \F)$, the set of equivalence classes of all irreducible admissible representations of $\GL(4, \F)$ is parametrized by the $\GL(4, \mathbb{C})$-conjugacy classes of all four-dimensional local Langlands parameters.~

Recall that the degree four symplectic group $\text{Sp}(4, \mathbb{C})$ is the complex dual group of $\SO(5, \F)$.~Let $\iota$ be the natural embedding of $\text{Sp}(4, \mathbb{C})$ into $\GL(4, \mathbb{C})$ and $\tau$ be an irreducible admissible representation of $\GL(4, \F)$ with $\varphi_{\tau}$ its local Langlands parameter under the local Langlands correspondence.~Then $\tau$ is a \textit{local Langlands functorial transfer} from $\SO(5, \F)$ to $\GL(4, \F)$ if $\varphi_{\tau}$ has image in a conjugate of $\text{Sp}(4, \mathbb{C})$.

\subsection{Local exterior square $L$-functions}\label{TwistedShalikaPeriod}
In this subsection, we introduce the function $\Lambda_{s_{0}}$ where $s_{0} \in \mathbb{C}$ and show how it may be used to compute local exterior square $L$-functions.~Next, we give a strategy for proving Theorem \ref{Main}.~

We denote the space of complex-valued, locally constant functions with compact support on $\F^{2}$ (space of \textit{Schwartz functions}) by $\mathcal{S}\left(\F^{2}\right)$.~Let $\tau$ be a supercuspidal representation of $\GL(4, \F)$.~For $s \in \mathbb{C}$, $W_{\tau} \in W\left(\tau, \psi_{\F}\right)$, and $\varphi \in \mathcal{S}\left(\F^{2}\right)$, we define 
\begin{align*}
    & J(s, W_{\tau}, \varphi, \psi_{\F}) \\
    &= \int\displaylimits_{\N(2, \, \F) \backslash \GL(2, \, \F)}\int\displaylimits_{\mathcal{B}(2, \, \F) \backslash \Mat(2 \times 2, \, \F)}W_{\tau}\left(\sigma_{4}\begin{pmatrix}
    I_{2} & X \\
    & I_{2}
\end{pmatrix}\begin{pmatrix}
    g & \\
    & g
\end{pmatrix}\right)\psi_{\F}\left(-\text{tr}X\right) |\det g|^{s}\varphi(\epsilon_{2}g)dXd^{\times}g
\end{align*}
where $\epsilon_{2} = (0, 1) \in \F^{2}$.~

The following theorem deals with the convergence of $J$ \cite[Section 7: Proposition 1, Section 9: Proposition 3]{Jacquet2}, \cite[Theorem 2.1]{YZ2}:
\begin{theorem}
There exists $r_{\tau, \, \wedge^{2}} \in \mathbb{R}$ such that for every $s \in \mathbb{C}$ with $\text{Re}(s) > r_{\tau, \, \wedge^{2}}$, $W_{\tau} \in W(\tau, \psi_{\F})$, and $\varphi \in \mathcal{S}\left(\F^{2}\right)$, the integral $J(s, W_{\tau}, \varphi, \psi_{\F})$ converges absolutely.
\end{theorem}
Using the preceding theorem, we may introduce the local exterior square $L$-function $L(s, \tau, \wedge^{2})$ \cite[Theorem 2.2]{YZ2}, \cite[Theorem 2.3]{Jo}, \cite[Lemma 3.1]{Cogdell2}:
\begin{theorem}
For fixed $W_{\tau} \in W(\tau, \psi_{\F})$ and $\varphi \in \mathcal{S}\left(\F^{2}\right)$, the map $s \mapsto J(s, W_{\tau}, \varphi, \psi_{\F})$ for $s \in \mathbb{C}$ with {\normalfont $\text{Re}(s) > r_{\tau, \, \wedge^{2}}$} defines an element of $\mathbb{C}\left(q_{\F}^{-s}\right)$ that is a rational function in the variable $q_{\F}^{-s}$ and, therefore, has a meromorphic continuation to the entire plane, which we continue to denote by $J(s, W_{\tau}, \varphi, \psi_{\F})$.~Furthermore, if we denote
\[
I_{\tau, \psi_{\F}} = \SP_{\mathbb{C}}\left\{J(s, W_{\tau}, \varphi, \psi_{\F}) : W_{\tau} \in W(\tau, \psi_{\F}), \; \varphi \in \mathcal{S}\left(\F^{2}\right)\right\},
\]
then there exists a unique $p(Z) \in \mathbb{C}[Z]$ such that $p(0) = 1$ and $I_{\tau, \psi_{\F}} = \frac{1}{p\left(q_{\F}^{-s}\right)}\mathbb{C}\left[q_{\F}^{-s}, q_{\F}^{s}\right]$. Moreover, $p(Z)$ does not depend on $\psi_{\F}$, and we denote $L(s, \tau, \wedge^{2}) = \frac{1}{p\left(q_{\F}^{-s}\right)}$.
\end{theorem}

For a more explicit description of $L(s,\tau,\wedge^2)$, let
$s_{0} \in \mathbb{C}$ satisfy
\[
\omega_\tau|\, \cdot \,|^{2s_0}\equiv 1
\]
where $\omega_{\tau}$ is the central character of $\tau$.~We introduce the function $\Lambda_{s_0}$ defined by \cite[Theorem 6.2 (i)]{Jo1}:
\begin{align}
    &\Lambda_{s_{0}}(W_{\tau}) \\ \notag
    & = \int\displaylimits_{\F^{\times}\N(2, \, \F) \backslash \GL(2, \, \F)}\int\displaylimits_{\mathcal{B}(2, \, \F) \backslash \Mat(2 \times 2, \, \F)}W_{\tau}\left(\sigma_{4}\begin{pmatrix}
    I_{2} & X \\
    & I_{2}
\end{pmatrix}\begin{pmatrix}
    g & \\
    & g
\end{pmatrix}\right)\psi_{\F}\left(-\text{tr}X\right) \, dX|\det g|^{s_{0}}d^{\times}g
\end{align}
where $W_{\tau} \in W(\tau, \psi_{\F})$.~Using the following result of Jo, one sees that for supercuspidal representations of $\GL(4, \F)$, their local exterior square $L$-functions are well understood \cite[Theorem 3.6 (i)]{Jo}, \cite[Theorem 6.2 (i)]{Jo1}.
\begin{theorem}\label{Jo}
    Let $\tau$ be a supercuspidal representation of $\GL(4, \F)$.~Then 
    \[
    L(s, \tau, \wedge^{2}) = \prod_{s_{0}} \frac{1}{1 - q_{\F}^{s_{0}-s}}
    \]
    with the product running over all $s_{0} \in \mathbb{C}$ such that $\omega_{\tau}|\, \cdot \,|^{2s_{0}} \equiv 1$ and $\Lambda_{s_{0}}(W_{\tau})$ is non-zero for some $W_{\tau} \in W(\tau, \psi_{\F})$.
\end{theorem}

From Jiang, Nien, and Qin, we have the following equivalent statements \cite[Theorem 5.5]{Jiang}:
\begin{theorem}\label{Jiang}
    Let $\tau$ be a supercuspidal representation of $\GL(4, \F)$.~The following are equivalent.~
    \begin{enumerate}
        \item $\tau$ is a local Langlands functorial transfer from $\SO(5, \F)$.
        \item $\tau$ has a non-zero local Shalika model.~
        \item The local exterior square $L$-function $L(s, \tau, \wedge^{2})$ has a pole at $s = 0$.
    \end{enumerate}
\end{theorem}
\noindent Statement (3) of Theorem \ref{Jiang} gives us a general strategy for proving Theorem \ref{Main}.~From Theorem \ref{Jo}, if $L(s, \tau, \wedge^{2})$ has a pole at $s = 0$, then $\omega_{\tau}$ is trivial.~Conversely, assume $\omega_{\tau}$ is trivial.~If we can exhibit $W_{\tau} \in W(\tau, \psi_{\F})$ such that $\Lambda_{0}(W_{\tau})$ is non-zero for middle and biquadratic supercuspidals, then Theorem \ref{Main} follows.~

\section{Construction of minimax supercuspidal representations induced from $[\mathfrak{A}, 1, 0, \beta]$}\label{ConstructionSupercuspidals}

In this section, we consider three classes of minimax supercuspidal representations that are induced from simple strata of the form $[\mathfrak{A}, 1, 0, \beta]$.~Furthermore, we give explicit Whittaker functions introduced in Subsection \ref{ExplicitWhittakerfunctions} associated with the supercuspidals we construct in subsections \ref{MiddleSupercuspidalRepresentations} and \ref{Biquadratic}.

\subsection{Simple supercuspidal representations}\label{SimpleSupercuspidalRepresentations}

In this subsection, we briefly recall the construction of \textit{simple supercuspidal representations} of $\GL(4, \F)$ via type theory (\cite{Knightly} and \cite[Example 2.18]{Ye}).~Let $\mathfrak{I}_{m}$ be the standard minimal hereditary $\mathcal{O}_{\F}$-order in $\Mat(4 \times 4, \F)$ with Jacobson radical $\mathfrak{P}_{m}$, and $v \in \mathcal{O}_{\F}^{\times}$.~We consider the simple stratum $\left[\mathfrak{I}_{m}, 1, 0, \beta_{v}\right]$ where
\[
\beta_{v} = \begin{pmatrix}
    & & & & \sfrac{1}{v\varpi_{\F}} \\
    1 & & & & \\
     & 1 & & & \\
      & & & 1 &
\end{pmatrix}
\]
is minimal over $\F$ and $\E_{v} = \F \left[ \beta_{v} \right]$ is a degree four totally ramified extension of $\F$.~

Associated with the simple stratum above are the open normal subgroups of $\mathfrak{I}_{m}^{\times}$:
\[
H^{1}\left(\beta_{v}, \mathfrak{I}_{m}\right) = J^{1}\left(\beta_{v}, \mathfrak{I}_{m}\right) = U^{1}\left(\mathfrak{I}_{m}\right) = I_{4} + \mathfrak{P}_{m}
\]
and the open compact subgroup $J_{v} := J\left(\beta_{v}, \mathfrak{I}_{m}\right) = \mathcal{O}_{\F}^{\times} \, U^{1}\left( \mathfrak{I}_{m}\right)$.~Let $\varphi$ be a quasi-character of $\mathcal{O}_{\F}^{\times}$ that is trivial on $1 + \mathcal{P}_{\F}$, i.e.~$\varphi$ is inflated from a quasi-character of $k_{\F}^{\times}$.~From the singleton $\mathcal{C}\left(\mathfrak{I}_{m}, \beta_{v}, \psi_{\F} \right) = \{\psi_{\beta_{v}}\}$, we may form a $\beta_{v}$-extension 
\[
\kappa_{(v,  \, \varphi)} \, \colon \, \mathcal{O}_{\F}^{\times}  \, U^{1}\left(\mathfrak{I}_{m}\right) \longrightarrow \mathbb{C}^{\times}
\]
of $\psi_{\beta_{v}}$ by
\[
\kappa_{(v,  \, \varphi)}(xy) := \varphi(x)\psi_{\beta_{v}}(y)
\]
where $x \in \mathcal{O}_{\F}^{\times}$ and $y \in U^{1}\left(\mathfrak{I}_{m}\right)$.

To form a maximal simple type with the data above, we take the pair $\left(J_{v},  \kappa_{(v,  \, \varphi)}\right)$.~Moreover, to form an extended maximal simple type $\left(\textbf{J}_{v},  \Lambda_{(v,  \, \varphi,  \, \zeta)}\right)$, let $\Lambda_{(v,  \, \varphi,  \, \zeta)}\left(\beta_{v}\right) = \zeta \in \mathbb{C}^{\times}$ so that $\Lambda_{(v,  \, \varphi,  \, \zeta)}(\varpi_{\F}) = \left(\zeta^{4}  \varphi(v)\right)^{-1}$.~Hence, we may now define a simple supercuspidal representation:
\begin{definition}
    \normalfont 
    Let $\left(\textbf{J}_{v},  \Lambda_{(v,  \, \varphi,  \, \zeta)}\right)$ be an extended maximal simple type as above.~A \textit{simple supercuspidal representation} of $\GL(4, \F)$ is a depth $\frac{1}{e\left(\mathfrak{I}_{m}\right)} = \frac{1}{4}$ minimax supercuspidal representation of the form 
    \[
    \pi_{(v,  \, \varphi,  \, \zeta)} = \ind_{\textbf{J}_{v}}^{\GL(4, \, \F)}\Lambda_{(v,  \, \varphi,  \, \zeta)}.
    \]
\end{definition}

\subsection{Middle supercuspidal representations}\label{MiddleSupercuspidalRepresentations}

In this subsection, we briefly recall the construction of \textit{middle supercuspidal representations} $\pi_{(f, \, \chi, \, \zeta)}$ of $\GL(4, \F)$ via type theory \cite[Subsection 3.1]{Luo}.~Next, we give an explicit Whittaker function $\mathcal{W}_{(f, \, \chi, \, \zeta)}$ in the Whittaker model $W\left(\pi_{(f, \, \chi, \, \zeta)}, \psi_{\F}\right)$ of $\pi_{(f, \, \chi, \, \zeta)}$ \cite[Subsection 3.2]{Luo}.~

Let $f = X^{2}-dX-c$ with $c$,  $d \in \mathcal{O}_{\F}$ such that $f\imod{\mathcal{P}_{\F}}$ is irreducible over $k_{\F}$.~Furthermore, let 
\[
g_{f} = \begin{pmatrix}
    1 & & & \\
    &  & 1 & \\
    & \sfrac{\varpi_{\F}}{c} &   & \\
    & & & \sfrac{\varpi_{\F}}{c}
\end{pmatrix} \in \GL(4, \F)
\]
and $\mathfrak{A}_{4}$ be the conjugate by $g_{f}$ of the standard hereditary $\mathcal{O}_{\F}$-order $\mathfrak{A}_{2}$ in $\Mat(4 \times 4, \F)$ corresponding to the partition $4 = 2+2$.~Explicitly, we have that
\begin{equation}\label{Order}
    \mathfrak{A}_{4} = \begin{pmatrix}
        \mathcal{O}_{\F} & \mathcal{O}_{\F}& \mathcal{P}_{\F}^{-1} & \mathcal{P}_{\F}^{-1}\\
        \mathcal{P}_{\F} & \mathcal{O}_{\F} & \mathcal{O}_{\F} & \mathcal{P}_{\F}^{-1}\\
        \mathcal{P}_{\F} & \mathcal{P}_{\F} & \mathcal{O}_{\F}& \mathcal{O}_{\F}\\
        \mathcal{P}_{\F}^{2} & \mathcal{P}_{\F} & \mathcal{P}_{\F}& \mathcal{O}_{\F}
    \end{pmatrix}
\end{equation}
with period $e\left(\mathfrak{A}_{4}\right) = 2$ and Jacobson radical 
\begin{equation}\label{Jacobson}
    \mathfrak{P}_{4} = \begin{pmatrix}
        \mathcal{P}_{\F} & \mathcal{O}_{\F}& \mathcal{O}_{\F} & \mathcal{P}_{\F}^{-1}\\
        \mathcal{P}_{\F} & \mathcal{P}_{\F} & \mathcal{O}_{\F} & \mathcal{O}_{\F}\\
        \mathcal{P}_{\F}^{2} & \mathcal{P}_{\F} & \mathcal{P}_{\F}& \mathcal{O}_{\F}\\
        \mathcal{P}_{\F}^{2} & \mathcal{P}_{\F}^{2} & \mathcal{P}_{\F}& \mathcal{P}_{\F}
    \end{pmatrix}.
\end{equation}

Next, let \begin{equation*}
    \beta_{f} := \begin{pmatrix}
    & & & \sfrac{c}{\varpi_{\F}^{2}} \\
    1 & & &  \\
     & 1 & & \sfrac{d}{\varpi_{\F}} \\
     & & 1&  
\end{pmatrix}.
\end{equation*}
The element $\beta_{f}$ gives us a degree four extension $\E_{f} = \F\left[\beta_{f}\right]$ of $\F$ with inertial degree and ramification index equal to two.~Furthermore, we have that $\nu_{\E_{f}}\left(\beta_{f}\right) = -1$ where $\nu_{\E_{f}}$ is the normalized discrete valuation on $\E_{f}$.~Our choices above give a simple stratum $\left[\mathfrak{A}_{4}, 1, 0, \beta_{f}\right]$ with $\beta_{f}$ minimal over $\F$.~This implies that the supercuspidal representations coming from our simple stratum will be minimax in the sense of \cite{Adrian}.~Furthermore, we let $\sigma_{f} := \beta_{f}^{2}\varpi_{\F}$ be a zero of $f$ and $\text{L}_{f} := \F\left[\sigma_{f}\right]$ the unique degree two unramified extension of $\F$.

Associated with our simple stratum is the open compact subgroup $J_{f} := J\left(\beta_{f}, \mathfrak{A}_{4}\right)$ with open normal subgroups $J^{1}:= J^{1}\left(\beta_{f}, \mathfrak{A}_{4}\right)$ and $H^{1}:= H^{1}\left(\beta_{f}, \mathfrak{A}_{4}\right)$ where
\[
J^{1} = H^{1} = U^{1}:= U^{1}\left(\mathfrak{A}_{4}\right) = I_{4} + \mathfrak{P}_{4}.
\]
Furthermore, $J_{f} = \mathcal{O}_{\text{L}_{f}}^{\times} \, U^{1}$ where if we fix the ordered $\mathcal{O}_{\F}$-basis $\{c,   \sigma_{f}\}$ of $\mathcal{O}_{\text{L}_{f}}$, a general element $x_{0}c+x_{1}\sigma_{f} \in \mathcal{O}_{\text{L}_{f}}^{\times}$ is of the form
\begin{equation}\label{Unramified}
\begin{pmatrix}
        x_{0}c & & \sfrac{x_{1}c}{\varpi_{\F}} & \\
        & x_{0}c & & \sfrac{x_{1}c}{\varpi_{\F}} \\
        x_{1}\varpi_{\F} & & x_{0}c+x_{1}d & \\
        & x_{1}\varpi_{\F} & & x_{0}c+x_{1}d
    \end{pmatrix}
\end{equation}
\noindent with $x_{i} \in \mathcal{O}_{\F}$.~

Let $\chi$ be a quasi-character of $\mathcal{O}_{\text{L}_{f}}^{\times}$ that is trivial on $1 + \mathcal{P}_{\text{L}_{f}}$, i.e.~$\chi$ is the inflation of a quasi-character of $k_{\text{L}_{f}}^{\times}$.~From $\chi$, we now recall simple types $\left(J_{f},  \kappa_{(f, \, \chi)}\right)$ associated with our simple stratum, where $\kappa_{(f, \, \chi)}$ is a $\beta_{f}$-extension of $J_{f}$.~Let $\psi_{\beta_{f}}$ be a quasi-character of $U^{1}$ of the form
\[
\psi_{\beta_{f}}(x) = \left(\psi_{\F} \circ \text{tr}_{\Mat(4 \times 4, \, \F)}\right)(\beta_{f}(x-1)).
\]
Since $\mathcal{C}\left(\mathfrak{A}_{4}, \beta_{f}, \psi_{\F} \right) = \{\psi_{\beta_{f}}\}$, we may form 
\[
\kappa_{(f, \, \chi)} \, \colon \, \mathcal{O}_{\text{L}_{f}}^{\times}  \, U^{1}\longrightarrow \mathbb{C}^{\times}
\]
by setting
\[
   \kappa_{(f, \, \chi)}(xy) :=  \chi(x)\psi_{\beta_{f}}(y),
\]
where $x \in \mathcal{O}_{\text{L}_{f}}^{\times}$ and $y \in U^{1}$.

To get an extended maximal simple type $\left(\textbf{J}_{f},  \Lambda_{(f,  \, \chi,  \, \zeta)}\right)$ from a maximal simple type defined above, we set $\Lambda_{(f,  \, \chi,  \, \zeta)}\left(\beta_{f}\right) = \zeta \in \mathbb{C}^{\times}$ so that $\Lambda_{(f,  \, \chi,  \, \zeta)}(\varpi_{\F}) = \zeta^{-2}\chi\left(\sigma_{f}\right)$.~This enables us to define a \textit{middle supercuspidal representation}.
\begin{definition}\label{MiddleSupercuspidal}
    \normalfont
    Let $\left(\textbf{J}_{f},  \Lambda_{(f,  \, \chi,  \, \zeta)}\right)$ be an extended maximal simple type as above.~A \textit{middle supercuspidal representation} of $\GL(4, \F)$ is a depth $\frac{1}{e\left(\mathfrak{A}_{4}\right)} = \frac{1}{2}$ minimax supercuspidal representation of the form 
    \[
    \pi_{(f,  \, \chi,  \, \zeta)} = \ind_{\textbf{J}_{f}}^{\GL(4, \, \F)}\Lambda_{(f,  \, \chi,  \, \zeta)}.
    \]
\end{definition}
\noindent Let $\mathcal{A}_{\text{middle}}^{(4)}$ denote the set of isomorphism classes of middle supercuspidal representations of $\GL(4, \F)$.~We recall the following bijection \cite[Proposition 3.2]{Luo}:
\begin{proposition}\label{Bijection}
    There exists a bijection between $\mathcal{A}_{\text{middle}}^{(4)}$ and the set of triples $(\bar{f},  \, \chi,  \, \zeta)$ with $\bar{f}$ a monic irreducible degree two polynomial over $k_{\F}$, $\chi$ a multiplicative quasi-character of $k_{\bar{f}} = \bigslant{k_{\F}[X]}{\left\langle \bar{f}(X) \right\rangle}$, and $\zeta \in \mathbb{C}^{\times}$.
\end{proposition}

Lastly, we recall the Whittaker function $\mathcal{W}_{\left(f, \, \chi, \, \zeta\right)} \, \colon \, \GL(4, \F) \longrightarrow \mathbb{C}$ where 
\begin{equation}
    \mathcal{W}_{\left(f, \, \chi, \, \zeta\right)}(g) = \psi_{4}(u)  \, \Lambda_{\left(f, \, \chi, \, \zeta\right)}(h)
\end{equation}
if $g = uh \in \N(4, \F)  \, \textbf{J}_{f}$ with $u \in \N(4, \F)$ and $h \in \textbf{J}_{f}$ and zero otherwise.~Furthermore, we have from \cite{Luo} that the support $\Su\left(\mathcal{W}_{(f,  \, \chi,  \, \zeta)}\right)$ of $\mathcal{W}_{(f,  \, \chi,  \, \zeta)}$ is contained inside the following disjoint union:
\begin{equation}\label{DisjointUnion}
\Su\left(\mathcal{W}_{(f,  \, \chi,  \, \zeta)}\right) \subset \displaystyle\bigsqcup_{k \in \mathbb{Z}}\N(4, \F)  \beta_{f}^{k}  J_{f}.
\end{equation}

\subsection{Biquadratic supercuspidal representations}\label{Biquadratic}

In this subsection, we construct depth-one minimax supercuspidal representations $\pi$ of $\GL(4, \F)$ and give an explicit associated Whittaker function as in Subsection \ref{ExplicitWhittakerfunctions}.~Let $f_{M} = X^{4}-bX^{2}-a$ with $a$,  $b \in \mathcal{O}_{\F}$ such that $\bar{f}_{M} = f_{M} \imod{\mathcal{P}_{\F}}$ is irreducible over $k_{\F}$ and let $t = \left(t_{i, j}\right) \in \GL(4,  \F)$ be the diagonal matrix defined by: $t_{i, i} = \varpi_{\F}^{i-1}$ for $1 \leq i \leq 4$.~We say $f_{M}$ is \textit{biquadratic} as it omits all odd power terms.~Furthermore, let $\mathfrak{A}_{M}$ be the conjugate by $t$ of the maximal hereditary $\mathcal{O}_{\F}$-order $\Mat(4\times 4, \mathcal{O}_{\F})$:
\begin{equation}\label{MaximalOrder}
    \mathfrak{A}_{M} = \begin{pmatrix}
        \mathcal{O}_{\F} & \mathcal{P}_{\F}^{-1} & \mathcal{P}_{\F}^{-2} & \mathcal{P}_{\F}^{-3} \\
        \mathcal{P}_{\F} & \mathcal{O}_{\F} & \mathcal{P}_{\F}^{-1} & \mathcal{P}_{\F}^{-2} \\
        \mathcal{P}_{\F}^{2} & \mathcal{P}_{\F} & \mathcal{O}_{\F} & \mathcal{P}_{\F}^{-1} \\
        \mathcal{P}_{\F}^{3} & \mathcal{P}_{\F}^{2} & \mathcal{P}_{\F} & \mathcal{O}_{\F}
    \end{pmatrix}.
\end{equation}

Next, let 
\begin{equation*}
    \beta_{f_{M}} := t\begin{pmatrix}
    & & & \sfrac{a}{\varpi_{\F}} \\
    \sfrac{1}{\varpi_{\F}} & & &  \\
     & \sfrac{1}{\varpi_{\F}} & & \sfrac{b}{\varpi_{\F}} \\
     & & \sfrac{1}{\varpi_{\F}} &  
\end{pmatrix}t^{-1} 
= \begin{pmatrix}
    & & & \sfrac{a}{\varpi_{\F}^{4}} \\
    1 & & &  \\
     & 1 & & \sfrac{b}{\varpi_{\F}^{2}} \\
     & & 1&  
\end{pmatrix}.
\end{equation*}
The element $\beta_{f_{M}}$ gives us a degree four unramified extension $\E_{f_{M}} = \F\left[\beta_{f_{M}}\right]$ of $\F$.~Furthermore, we have that $\nu_{\E_{f_{M}}}\left(\beta_{f_{M}}\right) = -1$ where $\nu_{\E_{f_{M}}}$ is the normalized discrete valuation on $\E_{f_{M}}$.~Let $\sigma_{f_{M}} := \beta_{f_{M}}\varpi_{\F}$ be a zero of $f_{M}$.~It is not hard to see that $\left[\mathfrak{A}_{M}, 1, 0, \beta_{f_{M}}\right]$ is a simple stratum with $\beta_{f_{M}}$ minimal over $\F$ \cite[Chapter 1]{BK}.~This implies that the supercuspidal representations coming from our simple stratum will be minimax in the sense of \cite{Adrian}.~
 
Associated with our simple stratum $\left[\mathfrak{A}_{M}, 1, 0, \beta_{f_{M}}\right]$ is the open compact subgroup $J_{f_{M}} := J\left(\beta_{f_{M}}, \mathfrak{A}_{M}\right)$ with open normal subgroups $J_{M}^{1}:= J_{M}^{1}\left(\beta_{f_{M}}, \mathfrak{A}_{M}\right)$ and $H_{M}^{1}:= H_{M}^{1}\left(\beta_{f_{M}}, \mathfrak{A}_{M}\right)$.~In this situation, we have that
\[
J_{M}^{1} = H_{M}^{1} = U_{M}^{1}:= U_{M}^{1}\left(\mathfrak{A}_{M}\right) = I_{4} + \mathfrak{P}_{M}
\]
where $\mathfrak{P}_{M} = \varpi_{\F}\mathfrak{A}_{M}$ is the Jacobson radical of $\mathfrak{A}_{M}$.~Furthermore, $J_{f_{M}} = \mathcal{O}_{\E_{f_{M}}}^{\times} \, U_{M}^{1}$ and if we fix the ordered $\mathcal{O}_{\F}$-basis $\left\{1, \sigma_{f_{M}}, \sigma_{f_{M}}^{2}, \sigma_{f_{M}}^{3}\right\}$ of $\mathcal{O}_{\E_{f_{M}}}$, a general element $x_{0}+x_{1}\sigma_{f_{M}}+x_{2}\sigma_{f_{M}}^{2}+x_{3}\sigma_{f_{M}}^{3} \in \mathcal{O}_{\E_{f_{M}}}^{\times}$ with $x_{i}\in \mathcal{O}_{\F}$ is of the form
\begin{equation}\label{Unramified2}
\begin{pmatrix}
       x_{0} & \sfrac{x_{3}a}{\varpi_{\F}} & \sfrac{x_{2}a}{\varpi_{\F}^{2}} & \sfrac{(x_{1}+x_{3}b)a}{\varpi_{\F}^{3}} \\
       x_{1}\varpi_{\F} & x_{0} &  \sfrac{x_{3}a}{\varpi_{\F}} & \sfrac{x_{2}a}{\varpi_{\F}^{2}} \\
       x_{2}\varpi_{\F}^{2} & (x_{1} + x_{3}b)\varpi_{\F} & x_{0}+x_{2}b & \sfrac{x_{1}b+x_{3}(a+b^{2})}{\varpi_{\F}} \\
       x_{3}\varpi_{\F}^{3} & x_{2}\varpi_{\F}^{2} & (x_{1} + x_{3}b)\varpi_{\F} & x_{0}+x_{2}b
    \end{pmatrix}.
\end{equation}

Next, we note that the set of simple characters $\mathcal{C}\left(\mathfrak{A}_{M}, \beta_{f_{M}}, \psi_{\F}\right)$ of $H_{M}^{1}$ is the singleton $\{\psi_{\beta_{f_{M}}}\}$ where
\[
\psi_{\beta_{f_{M}}}(x) = \left(\psi_{\F} \circ \text{tr}_{\Mat(4 \times 4, \, \F)}\right)(\beta_{f_{M}}(x-1))
\]
for $x \in U_{M}^{1}$.~From this simple character, we may define a $\beta_{f_{M}}$-\textit{extension} $\kappa_{(f_{M}, \, \chi_{M})}$ \cite[Chapter 5]{BK}.~In this situation, $\kappa_{(f_{M}, \, \chi_{M})}$ is a quasi-character of $J_{f_{M}}$ given by
\[
   \kappa_{(f_{M}, \, \chi_{M})}(xy) :=  \varphi_{M}(x)\psi_{\beta_{f_{M}}}(y)
\]
for $x \in \mathcal{O}_{\E_{f_{M}}}^{\times}$ and $y \in U_{M}^{1}$ where $\varphi_{M}$ is a quasi-character of $\mathcal{O}_{\E_{f_{M}}}^{\times}$ that is trivial on $1 + \mathcal{P}_{\E_{f_{M}}}^{2}$ and coincides with $\psi_{\beta_{f_{M}}}$ on $1 +  \mathcal{P}_{\E_{f_{M}}}$.~More explicitly, $\varphi_{M}$ is the inflation of a quasi-character $\bar{\varphi}_{M} \, \colon \, k_{\E_{f_{M}}}^{\times} \times k_{\E_{f_{M}}} \longrightarrow \mathbb{C}^{\times}$ of the form
\[
\bar{\varphi}_{M}(\bar{u}, \bar{z}) = \chi_{M}(\bar{u})\psi_{\beta_{f_{M}}}(1+\varpi_{\F}z)
\]
where $\chi_{M}$ is a quasi-character of $k_{\E_{f_{M}}}^{\times}$, $\bar{u} \in k_{\E_{f_{M}}}^{\times}$, and $z$ is any lift of $\bar{z} \in k_{\E_{f_{M}}}$.~This definition makes sense because $\psi_{\beta_{f_{M}}}$ is trivial on $1 + \mathcal{P}_{\E_{f_{M}}}^{2}$, hence, its restriction to $1 + \mathcal{P}_{\E_{f_{M}}}$ may be thought of as an additive quasi-character of $k_{\E_{f_{M}}}$.~

To see the inflation of $\bar{\varphi}_{M}$ to $\varphi_{M}$ more clearly, we first define a surjective homomorphism $\mathcal{O}_{\E_{f_{M}}}^{\times} \twoheadrightarrow k_{\E_{f_{M}}}^{\times} \times k_{\E_{f_{M}}}$ via
\[
x = u + uz\varpi_{\F} \mapsto \left(u \imod{1 + \mathcal{P}_{\E_{f_{M}}}}, \; z \imod{\mathcal{P}_{\E_{f_{M}}}} \right)
\]
where $u \in \mu_{\E_{f_{M}}}'$ and $z \in \mathcal{O}_{\E_{f_{M}}}$.~Hence, the kernel of this map is $1 + \mathcal{P}_{\E_{f_{M}}}^{2}$ and we have the isomorphism 
\[
\bigslant{\mathcal{O}_{\E_{f_{M}}}^{\times}}{1 + \mathcal{P}_{\E_{f_{M}}}^{2}} \cong k_{\E_{f_{M}}}^{\times} \times k_{\E_{f_{M}}}.
\]
There exist $q_{\F}^{4}-1$ distinct $\beta_{f_{M}}$-extensions and it follows trivially that $\kappa_{(f_{M}, \, \chi_{M})}$ satisfies the conditions for being a $\beta_{f_{M}}$-extension given in \cite[Chapter 5]{BK}.
 
We may now introduce simple types $\left(J_{f_{M}},  \lambda_{(f_{M}, \, \chi_{M})}\right)$ associated with our simple stratum.~In Bushnell--Kutzko notation, we remark that $\lambda_{(f_{M}, \, \chi_{M})} = \kappa_{(f_{M}, \, \chi_{M})} \otimes \sigma$ is a $\beta_{f_{M}}$-extension $\kappa_{(f_{M}, \, \chi_{M}')}$ as $\sigma$ is a cuspidal representation of $\bigslant{J_{f_{M}}}{J_{M}^{1}} \cong k_{\E_{f_{M}}}^{\times}$ which may be viewed as a quasi-character of $\mathcal{O}_{\E_{f_{M}}}^{\times}$ that is trivial on $1 + \mathcal{P}_{\E_{f_{M}}}$.~Moreover, we note that such a simple type is maximal as the period $e\left(\mathfrak{A}_{M}\right)$ of $\mathfrak{A}_{M}$ is equal to the ramification index one.

To get an extended maximal simple type $\left(\textbf{J}_{f_{M}},  \Lambda_{(f_{M},  \, \chi_{M},  \, \zeta)}\right)$ from a maximal simple type defined above, we set $\Lambda_{(f_{M},  \, \chi_{M},  \, \zeta)}\left(\beta_{f_{M}}\right) = \zeta \in \mathbb{C}^{\times}$ so that $\Lambda_{(f_{M},  \, \chi_{M},  \, \zeta)}(\varpi_{\F}) = \zeta^{-1}\chi_{M}\left(\sigma_{f_{M}}\right)$. From these choices, we may now define \textit{biquadratic supercuspidal representations}.
\begin{definition}\label{ONE}
    \normalfont
    We say $\pi_{(f_{M},  \, \chi_{M},  \, \zeta)}$ is a \textit{biquadratic supercuspidal representation} of $\GL(4,  \F)$ if 
    \[
    \pi_{(f_{M},  \, \chi_{M},  \, \zeta)} = \ind_{\textbf{J}_{f_{M}}}^{\GL(4, \,  \F)}\Lambda_{(f_{M},  \, \chi_{M},  \, \zeta)}
    \]
    where $\left(\textbf{J}_{f_{M}},  \Lambda_{(f_{M},  \, \chi_{M},  \, \zeta)}\right)$ is an extended maximal simple type of $\left(J_{f_{M}},  \kappa_{(f_{M}, \, \chi_{M})}\right)$.~Let $\omega_{(f_{M},  \, \chi_{M},  \, \zeta)}$ denote the central character of $\pi_{(f_{M},  \, \chi_{M},  \, \zeta)}$.~Furthermore, biquadratic supercuspidals have depth $\frac{1}{e\left(\mathfrak{A}_{M}\right)} = 1$.
\end{definition}

We may further classify these representations using the following proposition, which provides a parametrization for the set $\mathcal{A}_{\text{biquadratic}}$ of isomorphism classes of biquadratic supercuspidals of $\GL(4, \F)$.
\begin{proposition}\label{Bijection2}
    There exists a bijection between $\mathcal{A}_{\text{biquadratic}}$ and the set of triples $(\bar{f}_{M},  \, \chi_{M},  \, \zeta)$ with $\bar{f}_{M}$ a monic irreducible degree four biquadratic polynomial over $k_{\F}$, $\chi_{M}$ a quasi-character of $k_{\bar{f}_{M}}^{\times}$, where $ k_{\bar{f}_{M}}= \bigslant{k_{\F}[X]}{\left\langle \bar{f}_{M}(X) \right\rangle}$, and $\zeta \in \mathbb{C}^{\times}$.
\end{proposition}
\begin{proof}
    We first compute the minimal polynomial $\varphi_{\pi_{(f_{M},  \, \chi_{M},  \, \zeta)}}$ of $\left[\mathfrak{A}_{M}, 1, 0, \beta_{f_{M}}\right]$.~Let $\varphi_{Y_{\beta_{f_{M}}}}$ be the characteristic polynomial of $Y_{\beta_{f_{M}}} = \sigma_{f_{M}}$.~From Definition \ref{Characteristic}, we have that $\varphi_{\pi_{(f_{M},  \, \chi_{M},  \, \zeta)}}$ is equal to $\varphi_{Y_{\beta_{f_{M}}}}$ with coefficients reduced modulo $\mathcal{P}_{\F}$.~Hence, $\varphi_{\pi_{(f_{M},  \, \chi_{M},  \, \zeta)}}$ is a monic irreducible degree four biquadratic polynomial over $k_{\F}$.~The proof of Proposition \ref{Bijection2} now follows similarly from the proof of \cite[Proposition 3.2]{Luo}.~
\end{proof}

We end this subsection by giving an explicit Whittaker function $\mathcal{W}_{\left(f_{M}, \, \chi_{M}, \, \zeta\right)}$ in the Whittaker model $W \left(\pi_{\left(f_{M}, \, \chi_{M}, \, \zeta\right)}, \psi_{\F}\right)$ of $\pi_{(f_{M}, \, \chi_{M}, \, \zeta)}$ as in Subsection \ref{ExplicitWhittakerfunctions}.~We first note that
\begin{equation}\label{WhittakerBiQ}
    (J_{f_{M}} \cap \N(4, \F)) U^{1}_{M} = U^{1}_{M}.
\end{equation}
To see this, let $u \in J_{f_{M}} \cap \N(4, \F)$.~Then $u = xy$ for some $x \in \mathcal{O}_{\E_{f_{M}}}^{\times}$ and $y \in U^{1}_{M}$.~If $x^{-1}u \in U^{1}_{M}$, then comparing the first columns of $x^{-1}u$ and $y$ tells us that $x \in 1 + \mathcal{P}_{\E_{f_{M}}}$.~Hence, $u \in U^{1}_{M}$ and (\ref{WhittakerBiQ}) follows.~This implies $\Psi_{4} = \psi_{\beta_{f_{M}}}$.

Next, we compute the associated Bessel function $\mathcal{J}_{\left(f_{M}, \, \chi_{M}, \, \zeta\right)}$.~Let $\mathcal{N} = U^{2}_{M} := U^{2}\left(\mathfrak{A}_{M}\right)$. Then for $g \in \textbf{J}_{f_{M}}$, 
\begin{align*}
    \mathcal{J}_{\left(f_{M}, \, \chi_{M}, \, \zeta\right)}(g) &= \left(U^{1}_{M} : U^{2}_{M}  \right)^{-1} \displaystyle\sum_{h \in U^{1}_{M} \big/U^{2}_{M} } \notag \psi_{\beta_{f_{M}}}\left(h^{-1}\right)\Lambda_{\left(f_{M}, \, \chi_{M}, \, \zeta\right)}(gh)\\
    &= \Lambda_{\left(f_{M}, \, \chi_{M}, \, \zeta\right)}(g).~
\end{align*}
Hence, we may define $\mathcal{W}_{\left(f_{M}, \, \chi_{M}, \, \zeta\right)} \, \colon \, \GL(4, \F) \longrightarrow \mathbb{C}$ as 
\begin{align*}
    g \mapsto  \begin{cases} 
      \hfill \psi_{4}(u)  \Lambda_{\left(f_{M}, \, \chi_{M}, \, \zeta\right)}(h) &, \text{ if } g = uh \in \N(4, \F)  \textbf{J}_{f_{M}} \text{ with } u \in \N(4, \F), \; h \in \textbf{J}_{f_{M}} \\
      \hfill 0 & ,  \text{ otherwise}.
   \end{cases}
\end{align*}
Furthermore, since $\E_{f_{M}}$ is an unramified extension, $\varpi_{\F}$ is a uniformizer of $\mathcal{O}_{\E_{f_{M}}}$.~This implies that the support $\Su\left(\mathcal{W}_{(f_{M},  \, \chi_{M},  \, \zeta)}\right)$ of $\mathcal{W}_{(f_{M},  \, \chi_{M},  \, \zeta)}$ is contained inside the following disjoint union:
\begin{equation}\label{DisjointUnion2}
\Su\left(\mathcal{W}_{(f_{M},  \, \chi_{M},  \, \zeta)}\right) \subset \displaystyle\bigsqcup_{k \in \mathbb{Z}}\N(4, \F)  \varpi_{\F}^{k}  J_{f_{M}}.
\end{equation}

\section{Calculations}\label{ProofOfTheorem}

In this section, we give a proof of Theorem \ref{Main} following the strategy given at the end of Subsection \ref{TwistedShalikaPeriod}.~Let $\mathcal{N}_{2}^{-}$ denote the lower triangular nilpotent matrix subspace of $\Mat(2 \times 2, \F)$.~Then 
\begin{equation}\label{nilpotent}
    \mathcal{B}(2, \F) \backslash \Mat(2 \times 2, \F) \cong \mathcal{N}_{2}^{-} = \left\{n^{-}(x) = \begin{pmatrix}
        0 & 0 \\
        x & 0
    \end{pmatrix} \, : \, x \in \F \right\}.
\end{equation}
Furthermore, the Iwasawa decomposition of $\GL(2, \F)$ enables us to choose coset representatives of the form 
\begin{equation}\label{Iwasawa}
    \begin{pmatrix}
    \varpi_{\F}^{r} & \\
    & 1
\end{pmatrix}h
\end{equation}
with $r \in \mathbb{Z}$ and $h = (h_{i, j})  \in \GL(2, \mathcal{O}_{\F})$ for the quotient $\F^{\times}\N(2, \F)\backslash \GL(2, \F)$.~Let $g$ be of the form in (\ref{Iwasawa}) and $X = n^{-}(x)$ be of the form in (\ref{nilpotent}), then $\sigma_{4}\begin{pmatrix}
    I_{2} & n^{-}(x) \\
    & I_{2}
\end{pmatrix}\begin{pmatrix}
    g & \\
    & g
\end{pmatrix}\sigma_{4}$ is of the form 
\begin{equation}\label{Shalikaelement}
    \begin{pmatrix}
        h_{1, 1}\varpi_{\F}^{r} & & h_{1, 2}\varpi_{\F}^{r} & \\
         & h_{1, 1}\varpi_{\F}^{r} &  & h_{1, 2}\varpi_{\F}^{r} \\
         h_{2, 1} & h_{1, 1}x\varpi_{\F}^{r} & h_{2, 2} & h_{1, 2}x\varpi_{\F}^{r} \\
         & h_{2, 1} &  & h_{2, 2}
    \end{pmatrix}.
\end{equation}

Next, we mention a decomposition of the standard Iwahori subgroup $I$ of $\GL(2, \mathcal{O}_{\F})$ that will aid us in our computations.~Consider 
\[
I = I^{+}I_{0}I^{-}
\]
where $I^{+} = \N(2, \F) \cap I$, $I_{0} = A(2, \F) \cap I$, and $I^{-} = \N(2, \F)^{-} \cap I$ with $A(2, \F)$ denoting the standard diagonal torus of $\GL(2, \F)$ and $\N(2, \F)^{-}$ the unipotent radical of the Borel subgroup opposite to the standard Borel subgroup \cite[Section 3]{Chriss}.~This decomposition implies that if $h \in I$ in (\ref{Iwasawa}), then we may choose representatives of the form 
\begin{equation}\label{IwahoriReps}
     \begin{pmatrix}
    \varpi_{\F}^{r} & \\
    & 1
\end{pmatrix}\begin{pmatrix}
            h_{1, 1} & \\
            h_{2, 1} & h_{2, 2}
        \end{pmatrix}
\end{equation}
where $h_{1, 1}$, $h_{2, 2} \in \mathcal{O}_{\F}^{\times}$.

If $h \notin I$ in (\ref{Iwasawa}), then $h \in IsI$ from the Bruhat decomposition of $\GL(2, \mathcal{O}_{\F})$ where $s \in S_{2}$ is the longest Weyl element.~We may then use the decomposition of $I$ from the preceding paragraph to choose representatives of the form
\begin{equation}\label{IwahoriReps2}
\begin{pmatrix}
    \varpi_{\F}^{r} & \\
    & 1
\end{pmatrix}\begin{pmatrix}
    & h_{1, 2} \\
   h_{2, 1} & h_{2, 2}
\end{pmatrix}
\end{equation}
where $h_{1, 2}$, $h_{2, 1} \in \mathcal{O}_{\F}^{\times}$.~To see this, we first rewrite $IsI$ as $I^{+}sI$ since $sI_{0}I^{-}s \subset I$.~Next, we rewrite $I^{+}sI$ as $I^{+}(sI^{-}s)sI_{0}I^{+} = I^{+}sI_{0}I^{+}$ since $sI^{-}s \subset I^{+}$ \cite[Section 3]{Chriss}.~From this, our coset representative choices follow.~

\subsection{Calculation of $\Lambda_{0}\left(\pi_{(f, \, \chi, \, \zeta)}(\sigma_{4})\mathcal{W}_{(f, \, \chi, \, \zeta)}\right)$}

Let $\pi_{(f, \, \chi, \, \zeta)}$ be a middle supercuspidal representation of $\GL(4, \F)$.~Next, we consider the Whittaker function
\[
\pi_{(f, \, \chi, \, \zeta)}(\sigma_{4})\mathcal{W}_{(f, \, \chi, \, \zeta)} \in W(\pi_{(f, \, \chi, \, \zeta)}, \psi_{\F})
\]
where $\mathcal{W}_{(f, \, \chi, \, \zeta)}$ is the Whittaker function defined in Subsection \ref{MiddleSupercuspidalRepresentations}.~Our aim is to show that $\pi_{(f, \, \chi, \, \zeta)}(\sigma_{4})\mathcal{W}_{(f, \, \chi, \, \zeta)}$ satisfies the condition of our strategy given at the end of Subsection \ref{TwistedShalikaPeriod}.~We further recall from (\ref{DisjointUnion}) that $\Su\left(\mathcal{W}_{(f, \, \chi, \, \zeta)}\right)$ is contained inside the following: 
\[
\Su\left(\mathcal{W}_{(f, \, \chi, \, \zeta)}\right) \subset \displaystyle\bigsqcup_{k \in \mathbb{Z}}\N(4, \F)  \beta_{f}^{k}  J_{f}.
\]

For a general element $y$ of $J_{f}$, (\ref{Order}) and (\ref{Jacobson}) imply that we may write $y$ as
\begin{equation}\label{decomposition}
y = \begin{pmatrix}
        y_{1, 1} & y_{1, 2} & \sfrac{y_{1, 3}}{\varpi_{\F}} & \sfrac{y_{1, 4}}{\varpi_{\F}} \\ 
        y_{2, 1}\varpi_{\F} & y_{2, 2} & y_{2, 3} & \sfrac{y_{2, 4}}{\varpi_{\F}} \\ 
        y_{3, 1}\varpi_{\F} & y_{3, 2}\varpi_{\F} & y_{3, 3} & y_{3, 4}\\ 
        y_{4, 1}\varpi_{\F}^{2} & y_{4, 2}\varpi_{\F} & y_{4, 3}\varpi_{\F} & y_{4, 4}
    \end{pmatrix}   
    = (x_{0}c+x_{1}\sigma_{f})z
\end{equation}
where $x_{0}c+x_{1}\sigma_{f} \in \mathcal{O}_{\text{L}_{f}}^{\times}$ is of the form in (\ref{Unramified}) and $z \in U^{1}$ is of the form
\[
z = \begin{pmatrix}
        1+z_{1, 1}\varpi_{\F} & z_{1, 2} & z_{1, 3} & \sfrac{z_{1, 4}}{\varpi_{\F}} \\ 
        z_{2, 1}\varpi_{\F} & 1+z_{2, 2}\varpi_{\F} & z_{2, 3} & z_{2, 4} \\ 
        z_{3, 1}\varpi_{\F}^{2} & z_{3, 2}\varpi_{\F} & 1+z_{3, 3}\varpi_{\F} & z_{3, 4}\\ 
        z_{4, 1}\varpi_{\F}^{2} & z_{4, 2}\varpi_{\F}^{2} & z_{4, 3}\varpi_{\F} & 1+z_{4, 4}\varpi_{\F} 
    \end{pmatrix}  \in U^{1}
\]
with $y_{i, j}$, $z_{i, j} \in \mathcal{O}_{\F}$.~From (\ref{decomposition}), we note the following:
\begin{equation}\label{DETYM}
    \det(y) \equiv (y_{1, 3}y_{3, 1}-y_{1, 1}y_{3, 3})(y_{2, 4}y_{4, 2}-y_{2, 2}y_{4, 4}) \imod{\mathcal{P}_{\F}}.
\end{equation}
\begin{lemma}\label{Lemma}
    Let $g$ be of the form in (\ref{Iwasawa}) and $n^{-}(x)$ be of the form in (\ref{nilpotent}).~Next, set $\alpha := \sigma_{4}\begin{pmatrix}
    I_{2} & n^{-}(x) \\
    & I_{2}
\end{pmatrix}\begin{pmatrix}
    g & \\
    & g
\end{pmatrix}\sigma_{4}$.~If $\alpha \in \Su\left(\mathcal{W}_{(f, \, \chi, \, \zeta)}\right)$, then 
\[
\alpha \in \N(4, \F)J_{f} \, \sqcup \, \N(4, \F)\beta_{f}^{2}J_{f}.
\]
\end{lemma}
\begin{proof}
    Let $g$ and $n^{-}(x)$ be as in the hypothesis, then $\alpha$ is of the form in (\ref{Shalikaelement}).~Suppose $k <0$.~Then $\alpha = u\beta_{f}^{k}y$ for some $u \in \N(4, \F)$ and $y \in J_{f}$.~Comparing the $(4, 2)$ and $(4, 4)$-entries of $\alpha$ and $u\beta_{f}^{k}y$, we have that $h_{2, 1}$, $h_{2, 2} \in \mathcal{P}_{\F}$ which implies $h \notin \GL(2, \mathcal{O}_{\F})$.~This is a contradiction and $k \geq 0$.~

Suppose $k = 1$.~Comparing the $(4, 1)$ and $(4, 3)$-entries of $\alpha$ and $u\beta_{f}y$, we have that $y_{3, 1} = y_{3, 3} = 0$ which implies $\det(y) \in \mathcal{P}_{\F}$ from (\ref{DETYM}).~This implies $k \neq 1$.~For $k \neq 1$ and $k \equiv 1, 3 \imod{4}$, we have from comparing the $(4, 1)$ and $(4, 3)$-entries of $\alpha$ and $u\beta_{f}^{k}y$ the equations
\begin{align*}
(\beta_{f}^{k})_{4, 1}y_{1, 1} + (\beta_{f}^{k})_{4, 3}y_{3, 1}\varpi_{\F} &= 0 \\
(\beta_{f}^{k})_{4, 1}y_{1, 3}\varpi_{\F}^{-1} + (\beta_{f}^{k})_{4, 3}y_{3, 3} &= 0.
\end{align*}
Since $k \equiv 1, 3 \imod{4}$, at least one of $(\beta_{f}^{k})_{4, 1}$ or $(\beta_{f}^{k})_{4, 3}$ is non-zero.~Hence, the equations above imply $k \equiv 0, 2 \imod{4}$ as $\det(y) \in \mathcal{P}_{\F}$ otherwise from (\ref{DETYM}).

Suppose $k > 2$.~For $i \in \{2, 4\}$, we have that $\left(\beta_{f}^{k}\right)_{4, i} = C_{i}\varpi_{\F}^{-j}$ for some $j \geq 1$ and $C_{i} \in \mathcal{O}_{\F}$ where at least one of the $C_{i}$ is non-zero.~Comparing the $(4, 2)$ and $(4, 4)$-entries of $\alpha$ and $u\beta_{f}^{k}y$, we have that $(u\beta_{f}^{k}y)_{4, 2}$, $(u\beta_{f}^{k}y)_{4, 4} \in \mathcal{O}_{\F}$.~This implies $C_{2}y_{2, 2} + C_{4}y_{4, 2}$, $C_{2}y_{2, 4} + C_{4}y_{4, 4} \in \mathcal{P}_{\F}$.~Hence, this proves the lemma as $\det(y) \in \mathcal{P}_{\F}$ for $k > 2$ from (\ref{DETYM}).~
\end{proof}

\begin{proposition}\label{middle}
    Let $\pi_{(f, \, \chi, \, \zeta)}$ be a middle supercuspidal representation of $\GL(4, \F)$.~Then $L\left(s, \pi_{(f, \, \chi, \, \zeta)}, \wedge^{2}\right)$ has a pole at $s = 0$ if and only if $\omega_{(f, \, \chi, \, \zeta)}$ is trivial.
\end{proposition}

\begin{proof}
    If $L\left(s, \pi_{(f, \, \chi, \, \zeta)}, \wedge^{2}\right)$ has a pole at $s = 0$, then $\omega_{(f, \, \chi, \, \zeta)}$ is necessarily trivial from Theorem \ref{Jo}.~To prove the converse, suppose $\omega_{(f, \, \chi, \, \zeta)}$ is trivial.~We claim that $\Lambda_{0}$ evaluated at $\pi_{(f, \, \chi, \, \zeta)}(\sigma_{4})\mathcal{W}_{(f, \, \chi, \, \zeta)}$ is non-zero.
    
    Using our coset representative choices in (\ref{nilpotent}) and (\ref{Iwasawa}), we have that 
    \begin{align*}
        & \Lambda_{0}\left(\pi_{(f, \, \chi, \, \zeta)}(\sigma_{4})\mathcal{W}_{(f, \, \chi, \, \zeta)}\right) \\
        &= \int\displaylimits_{\F^{\times}\N(2, \, \F) \backslash \GL(2, \, \F)}\int\displaylimits_{\F}\mathcal{W}_{(f, \, \chi, \, \zeta)}\left(\sigma_{4}\begin{pmatrix}
    I_{2} & n^{-}(x) \\
    & I_{2}
\end{pmatrix}\begin{pmatrix}
    g & \\
    & g
\end{pmatrix}\sigma_{4}\right)\, dx \, d^{\times}g.
    \end{align*}
Furthermore, Lemma \ref{Lemma} tells us that we need only consider the double cosets 
\[
\F^{\times}\N(2, \F)\begin{pmatrix}
    \varpi_{\F}^{-2k} & \\
    & 1
\end{pmatrix}\GL(2, \mathcal{O}_{\F})
\]
with $0 \leq k \leq 1$ for the overall value of $\Lambda_{0}\left(\pi_{(f, \, \chi, \, \zeta)}(\sigma_{4})\mathcal{W}_{(f, \, \chi, \, \zeta)}\right)$.~Furthermore, for each $0 \leq k \leq 1$, if an element $\alpha$ of the form in (\ref{Shalikaelement}) with respect to $\F^{\times}\N(2, \F)\begin{pmatrix}
    \varpi_{\F}^{-2k} & \\
    & 1
\end{pmatrix}\GL(2, \mathcal{O}_{\F})$ is contained inside $\Su\left(\mathcal{W}_{(f, \, \chi, \, \zeta)}\right)$, then $\alpha$ is contained inside $\N(4, \F)\beta_{f}^{2k}J_{f}$ as $\det(\alpha) \in \mathcal{P}_{\F}^{-4k}$.~Hence, we may consider two separate integrals $T_{0}$ and $T_{1}$ such that 
\begin{equation}\label{TwoIntegrals}
\Lambda_{0}\left(\pi_{(f, \, \chi, \, \zeta)}(\sigma_{4})\mathcal{W}_{(f, \, \chi, \, \zeta)}\right) =T_{0} + T_{1}
\end{equation}
where 
\[
T_{k} = \int\displaylimits_{\F^{\times}\N(2, \, \F) \backslash \F^{\times}\N(2, \, \F)\begin{pmatrix}
    \varpi_{\F}^{-2k} & \\
    & 1
\end{pmatrix}\GL(2, \, \mathcal{O}_{\F})}\int\displaylimits_{\F}\mathcal{W}_{(f, \, \chi, \, \zeta)}\left(\sigma_{4}\begin{pmatrix}
    I_{2} & n^{-}(x) \\
    & I_{2}
\end{pmatrix}\begin{pmatrix}
    g & \\
    & g
\end{pmatrix}\sigma_{4}\right)\, dx \, d^{\times}g
\]
for $0 \leq k \leq 1$.

We first compute $T_{0}$.~If an element $\alpha$ of the form in (\ref{Shalikaelement}) is contained inside $\N(4, \F)J_{f}$, then $\alpha = uy$ for some $u = (u_{i, j}) \in \N(4, \F)$ and $y=(x_{0}c+x_{1}\sigma_{f})z \in J_{f}$.~Since $\det(\alpha) = \det(y) \in \mathcal{O}_{\F}^{\times}$, we may choose representatives of the form in (\ref{Iwasawa}) such that $r=0$.~Comparing the $(3, 1)$ and $(4, 2)$-entries of $\alpha$ and $uy$, we have that $h_{2, 1} \in (x_{1}+\mathcal{P}_{\F})\varpi_{\F}$ from (\ref{decomposition}); this implies $h \in I$.~Using (\ref{IwahoriReps}), we may choose representatives of the form 
        \[
        \begin{pmatrix}
            h_{1, 1} & \\
            h_{2, 1} & h_{2, 2}
        \end{pmatrix} \in I
        \]
        where $h_{1, 1}$, $h_{2, 2} \in \mathcal{O}_{\F}^{\times}$.~Comparing the $(3, 3)$ and $(4, 4)$-entries of $\alpha$ and $uy$, we see that $h_{2, 2}\in x_{0}c + x_{1}d +\mathcal{P}_{\F}$ from (\ref{decomposition}).

       To see where $x$ lies, we first suppose that $h_{2, 1}$ is zero.~Comparing the $(3, 2)$-entries of $\alpha$ and $uy$, we have that $x \in \mathcal{P}_{\F}$ as $h_{1, 1} \in \mathcal{O}_{\F}^{\times}$.~Suppose $h_{2, 1}$ is non-zero.~Comparing the $(3, 2)$ and $(3, 4)$-entries of $\alpha$ and $uy$, we have
        \begin{equation}\label{T0M1}
            u_{3, 4} = -\frac{y_{3, 4}}{h_{2, 2}} = \frac{h_{1, 1} x - y_{3, 2}\varpi_{\F}}{h_{2, 1}}.
        \end{equation}
    Hence, equation (\ref{T0M1}) implies
        \[
            x = \frac{h_{2, 2}y_{3, 2}\varpi_{\F}-h_{2, 1}y_{3, 4}}{h_{1, 1} h_{2, 2} } \in \mathcal{P}_{\F}.
        \]
        
To see where $h_{1, 1}$ lies, we first compare the $(2, 3)$-entries of $\alpha$ and $uy$ and see that $u_{2, 3} \in \mathcal{O}_{\F}$ as $h_{2, 2} \in \mathcal{O}_{\F}^{\times}$.~From this, we compare the $(2, 2)$ and $(2, 4)$-entries of $\alpha$ and $uy$ to get the equation
        \begin{equation}\label{T0M2}
            u_{2, 4} = \frac{h_{1, 1}-y_{2, 2}-u_{2, 3}y_{3, 2}\varpi_{\F}}{h_{2, 1}} = -\frac{y_{2, 4}+u_{2, 3}y_{3, 4}\varpi_{\F}}{h_{2, 2}\varpi_{\F}}.
        \end{equation}
        Equations (\ref{T0M2}) and (\ref{decomposition}) imply 
        \[
        h_{1, 1} \in \frac{(x_{0}^{2}c+x_{0}x_{1}d-x_{1}^{2})c}{x_{0}c + x_{1}d} +\mathcal{P}_{\F}.
        \]
        Hence, $\alpha$ lies inside
    \[
    \begin{pmatrix}
        \frac{(x_{0}^{2}c+x_{0}x_{1}d-x_{1}^{2})c}{x_{0}c + x_{1}d}+\mathcal{P}_{\F} & &  & \\
        & \frac{(x_{0}^{2}c+x_{0}x_{1}d-x_{1}^{2})c}{x_{0}c + x_{1}d}+\mathcal{P}_{\F} & &   \\
        (x_{1}+\mathcal{P}_{\F})\varpi_{\F} & \mathcal{P}_{\F} & x_{0}c + x_{1}d + \mathcal{P}_{\F} &  \\
        & (x_{1}+\mathcal{P}_{\F})\varpi_{\F} & & x_{0}c + x_{1}d + \mathcal{P}_{\F}
    \end{pmatrix}.
    \]
    
    Our computations above imply that we may decompose $\alpha$ as 
\[
\alpha = u(x_{0}c+x_{1}\sigma_{f})z,
\]
where the entries of $u$ satisfy
\[
    u_{1, 3} = u_{2, 4} = -\frac{x_{1}c}{h_{2, 2}\varpi_{\F}}, \; \;
    u_{1, 2} = u_{1, 4} = u_{2, 3} = u_{3, 4} = 0
\]
 and the entries of $z$ satisfy
\begin{align*}
 & 1 + z_{1, 1}\varpi_{\F}=1 + z_{2, 2}\varpi_{\F} = \frac{h_{2, 1}x_{1}c(x_{0}c+x_{1}d-h_{2, 2})+h_{1, 1}h_{2, 2}(x_{0}c+x_{1}d)\varpi_{\F}}{h_{2, 2}c(x_{0}^{2}c+x_{0}x_{1}d-x_{1}^{2})\varpi_{\F}} \\
 & z_{1, 2} = \frac{h_{1, 1}xx_{1}(x_{0}c+x_{1}d-h_{2, 2})}{h_{2, 2}(x_{0}^{2}c+x_{0}x_{1}d-x_{1}^{2})\varpi_{\F}}, \; \;
     z_{1, 3} = z_{2, 4} = \frac{x_{1}(x_{0}c+x_{1}d-h_{2, 2})}{(x_{0}^{2}c+x_{0}x_{1}d-x_{1}^{2})\varpi_{\F}} \\
      & z_{3, 1}\varpi_{\F}^{2}=z_{4, 2}\varpi_{\F}^{2} = \frac{h_{2, 1}c(h_{2, 2}x_{0}-x_{1}^{2})-h_{1, 1}h_{2, 2}x_{1}\varpi_{\F}}{h_{2, 2}c(x_{0}^{2}c+x_{0}x_{1}d-x_{1}^{2})}, \; \;
      z_{3, 2}\varpi_{\F}  = \frac{h_{1, 1}x(h_{2, 2}x_{0}-x_{1}^{2})}{h_{2, 2}(x_{0}^{2}c+x_{0}x_{1}d-x_{1}^{2})} \\
      & 1+ z_{3, 3}\varpi_{\F}= 1+z_{4, 4}\varpi_{\F} = \frac{h_{2, 2}x_{0}-x_{1}^{2}}{x_{0}^{2}c+x_{0}x_{1}d-x_{1}^{2}}
\end{align*}
with $z_{i, j}=0$ for all other $i$ and $j$.~For the entries of $z$, we note that $x_{0}^{2}c+x_{0}x_{1}d-x_{1}^{2} \in \mathcal{O}_{\F}^{\times}$ as $\det(x_{0}c+x_{1}\sigma_{f}) = ((x_{0}^{2}c+x_{0}x_{1}d-x_{1}^{2})c)^{2}  \in \mathcal{O}_{\F}^{\times}$.

By definition of $\mathcal{W}_{(f, \, \chi, \, \zeta)}$, we have that 
\begin{equation}\label{WhittakerMiddle1}
\mathcal{W}_{(f, \, \chi, \, \zeta)}(\alpha) =\chi(x_{0}c+x_{1}\sigma_{f})
\end{equation}
as $\psi_{4}(u) = \psi_{\beta_{f}}(z) = 1$.~Hence, equation (\ref{WhittakerMiddle1}) implies that $T_{0}$ is zero if $\chi$ is non-trivial as we are integrating over a compact subgroup $I$ and $\mathcal{P}_{\F}$.~If $\chi$ is trivial, then $T_{0}$ is equal to a positive real volume.

Next, we compute $T_{1}$.~If an element $\alpha$ of the form in (\ref{Shalikaelement}) is contained inside $\N(4, \F)\beta_{f}^{2}J_{f}$, then $\alpha = u\beta_{f}^{2}y$ for some $u = (u_{i, j}) \in \N(4, \F)$ and $y=(x_{0}c+x_{1}\sigma_{f})z \in J_{f}$.~Since $\det(\alpha) = (c\varpi_{\F}^{-2})^{2}\det(y)$, we may choose representatives of the form in (\ref{Iwasawa}) such that $r = -2$.~Comparing the $(3, 3)$ and $(4, 4)$-entries of $\alpha$ and $u\beta_{f}^{2}y$, we have that $y_{1, 3} \equiv -dy_{3, 3} \imod{\mathcal{P}_{\F}}$ and $y_{2, 4} \equiv -dy_{4, 4} \imod{\mathcal{P}_{\F}}$.~This implies $\det(y) \equiv (x_{0}c + x_{1}d)^{4}\imod{\mathcal{P}_{\F}}$ from (\ref{DETYM}) and (\ref{decomposition}).~Hence, $x_{0}c + x_{1}d \in \mathcal{O}_{\F}^{\times}$.~Comparing the $(3, 1)$ and $(4, 2)$-entries of $\alpha$ and $u\beta_{f}^{2}y$, we have that $h_{2, 1} \in x_{0}c + x_{1}d + \mathcal{P}_{\F} \subset \mathcal{O}_{\F}^{\times}$ from (\ref{decomposition}); this implies $h \in IsI$.~    
Using (\ref{IwahoriReps2}), we may choose representatives of the form 
\[
\begin{pmatrix}
    \varpi_{\F}^{-2} & \\
    & 1
\end{pmatrix}\begin{pmatrix}
    & h_{1, 2} \\
   h_{2, 1} & h_{2, 2}
\end{pmatrix}
\]
where $h_{1, 2}$, $h_{2, 1} \in \mathcal{O}_{\F}^{\times}$.~

To see where $x$ lies, we compare the $(3, 2)$ and $(3, 4)$-entries of $\alpha$ and $u\beta_{f}^{2}y$ and get the equation
\begin{equation}\label{T0M3}
    u_{3, 4} = \frac{h_{1, 2}x-y_{1, 4}\varpi_{\F}-dy_{3, 4}\varpi_{\F}}{h_{2, 2}\varpi_{\F}^{2}} = -\frac{y_{1, 2}+dy_{3, 2}}{h_{2, 1}}.
\end{equation}
Hence, equation (\ref{T0M3}) implies
        \[
            x = \frac{(h_{2, 1}(y_{1, 4}+dy_{3, 4})-h_{2, 2}(y_{1, 2}+dy_{3, 2})\varpi_{\F})\varpi_{\F}}{h_{1, 2} h_{2, 1}} \in \mathcal{P}_{\F}.
        \]

To see where $h_{1, 2}$ lies, we first compare the $(2, 1)$-entries of $\alpha$ and $u\beta_{f}^{2}y$.~In doing so, we see that $u_{2, 3} \in \mathcal{O}_{\F}$ as $h_{2, 1} \in \mathcal{O}_{\F}^{\times}$.~From this, we compare the $(2, 2)$ and $(2, 4)$-entries of $\alpha$ and $u\beta_{f}^{2}y$ to get the equation 
\begin{equation}\label{T0M4}
    u_{2, 4} = -\frac{cy_{4, 2}+u_{2, 3}(y_{1, 2}+dy_{3, 2})\varpi_{\F}}{h_{2, 1}\varpi_{\F}} = \frac{h_{1, 2}-(cy_{4, 4}+u_{2, 3}(y_{1, 4}+dy_{3, 4})\varpi_{\F})}{h_{2, 2}\varpi_{\F}^{2}}.
\end{equation}
Equations (\ref{T0M4}) and (\ref{decomposition}) imply $h_{1, 2} \in (x_{0}c + x_{1}d)c +\mathcal{P}_{\F}$.~Hence, $\alpha$ lies inside
\[
\begin{pmatrix}
    & & ((x_{0}c + x_{1}d)c+ \mathcal{P}_{\F})\varpi_{\F}^{-2}& \\
    & & & ((x_{0}c + x_{1}d)c + \mathcal{P}_{\F})\varpi_{\F}^{-2} \\
   x_{0}c + x_{1}d + \mathcal{P}_{\F} & & \mathcal{O}_{\F} & \mathcal{P}_{\F}^{-1} \\
    & x_{0}c + x_{1}d + \mathcal{P}_{\F} & & \mathcal{O}_{\F}
\end{pmatrix}.
\]

Returning to the $(4, 4)$-entries of $\alpha$ and $u\beta_{f}^{2}y$, we see that $x_{1}c \equiv -(x_{0}c + x_{1}d)d \imod{\mathcal{P}_{\F}}$ since $(u\beta_{f}^{2}y)_{4, 4}\in \mathcal{O}_{\F}$.~This implies 
\begin{equation}\label{EquivalentOne}
    x_{0} \equiv -\left(\frac{c+d^{2}}{cd}\right)x_{1} \imod{\mathcal{P}_{\F}}.
\end{equation}
Our computations above imply that we may decompose $\alpha$ as
\[
\alpha = u\beta_{f}^{2}x_{1}\left(-\frac{c+d^{2}}{d}+\sigma_{f}\right)z,
\]
where the entries of $u$ satisfy
\[
u_{1, 3} = u_{2, 4} = -\frac{c(h_{2, 1}d+x_{1}(c+d^{2}))}{h_{2, 1}d^{2}\varpi_{\F}}, \; \; u_{1, 2} = u_{1, 4} = u_{2, 3} = u_{3, 4} = 0
\]
and the entries of $z$ satisfy
\begin{align*}
  & 1+ z_{1, 1}\varpi_{\F}= 1+z_{2, 2}\varpi_{\F} = -\frac{h_{2, 1}d}{x_{1}c}, \; \;
  z_{1, 3} = z_{2, 4} = -\frac{h_{2, 2}d}{x_{1}c}, \; \;
z_{1, 4}\varpi_{\F}^{-1} = -\frac{h_{1, 2}xd}{x_{1}c\varpi_{\F}^{2}} \\
& z_{3, 1}\varpi_{\F}^{2}=z_{4, 2}\varpi_{\F}^{2} = -\frac{(c+d^{2})(h_{2, 1}d+x_{1}c)\varpi_{\F}}{x_{1}c^{2}d}, \; \; z_{3, 4} =  -\frac{h_{1, 2}x(c+d^{2})(h_{2, 1}d+x_{1}c)}{h_{2, 1}x_{1}c^{2}d\varpi_{\F}} \\
     & 1+ z_{3, 3}\varpi_{\F} = 1+z_{4, 4}\varpi_{\F} =-\frac{h_{1, 2}h_{2, 1}d^{2}+h_{2, 2}(c+d^{2})(h_{2, 1}d+x_{1}c)\varpi_{\F}}{h_{2, 1}x_{1}c^{2}d}
\end{align*}
with $z_{i, j} = 0$ for all other $i$ and $j$.~

By definition of $\mathcal{W}_{(f, \, \chi, \, \zeta)}$, we have that 
\begin{equation}\label{WhittakerMiddle2}
    \mathcal{W}_{(f, \, \chi, \, \zeta)}(\alpha) =\zeta^{2} \, \chi\left(-\frac{c+d^{2}}{d}+\sigma_{f}\right)
\end{equation}
as $\omega_{(f, \, \chi, \, \zeta)}$ is trivial and $\psi_{4}(u) = \psi_{\beta_{f}}(z) = 1$.~Since $\sigma_{f}^{2} = d\sigma_{f}+c$, (\ref{WhittakerMiddle2}) simplifies to
\begin{align}\label{Final}
    \zeta^{2} \, \chi\left(-\frac{c+d^{2}}{d}+\sigma_{f}\right) &= \chi\left(\left(-\frac{c+d^{2}}{d}\right)\sigma_{f} +\sigma_{f}^{2}\right) \\
    \notag &= \chi\left(\frac{c(d - \sigma_{f})}{d}\right).
\end{align}

Equation (\ref{Final}) implies that $T_{1}$ is equal to $\chi\left(\frac{c(d - \sigma_{f})}{d}\right)\text{V}$ for some positive real volume $\text{V}$ as we are integrating over a compact subgroup $IsI$ and $\mathcal{P}_{\F}$.~Our calculations show that in either situation with $\chi$ trivial or non-trivial,   $\Lambda_{0}$ evaluated at $\pi_{(f, \, \chi, \, \zeta)}(\sigma_{4})\mathcal{W}_{(f, \, \chi, \, \zeta)}$ is non-zero from (\ref{TwoIntegrals}).~The proof of the converse now follows from Theorem \ref{Jo}.
\end{proof}

\subsection{Calculation of $\Lambda_{0}\left(\pi_{(f_{M}, \, \chi_{M}, \, \zeta)}(\sigma_{4})\mathcal{W}_{(f_{M}, \, \chi_{M}, \, \zeta)}\right)$}

Let $\pi_{(f_{M}, \, \chi_{M}, \, \zeta)}$ be a biquadratic supercuspidal representation of $\GL(4, \F)$.~Next, we consider the Whittaker function
\[
\pi_{(f_{M}, \, \chi_{M}, \, \zeta)}(\sigma_{4})\mathcal{W}_{(f_{M}, \, \chi_{M}, \, \zeta)} \in W(\pi_{(f_{M}, \, \chi_{M}, \, \zeta)}, \psi_{\F})
\]
where $\mathcal{W}_{(f_{M}, \, \chi_{M}, \, \zeta)}$ is the Whittaker function defined in Subsection \ref{Biquadratic}.~Our aim is to show that $\pi_{(f_{M}, \, \chi_{M}, \, \zeta)}(\sigma_{4})\mathcal{W}_{(f_{M}, \, \chi_{M}, \, \zeta)}$ satisfies the condition of our strategy given at the end of Subsection \ref{TwistedShalikaPeriod}.~We further recall from (\ref{DisjointUnion2}) that $\Su\left(\mathcal{W}_{(f_{M}, \, \chi_{M}, \, \zeta)}\right)$ is contained inside the following:
\begin{equation*}
    \Su\left(\mathcal{W}_{(f_{M}, \, \chi_{M}, \, \zeta)}\right) \subset \displaystyle\bigsqcup_{k \in \mathbb{Z}}\N(4, \F)  \varpi_{\F}^{k}  J_{f_{M}}.
\end{equation*}

For a general element $y$ of $J_{f_{M}}$, (\ref{MaximalOrder}) and (\ref{Unramified2}) imply that we may write a general element $y$ as
\begin{equation}\label{decomposition2}
y = \begin{pmatrix}
        y_{1, 1} & \sfrac{y_{1, 2}}{\varpi_{\F}} & \sfrac{y_{1, 3}}{\varpi_{\F}^{2}} & \sfrac{y_{1, 4}}{\varpi_{\F}^{3}} \\ 
        y_{2, 1}\varpi_{\F} & y_{2, 2} & \sfrac{y_{2, 3}}{\varpi_{\F}} & \sfrac{y_{2, 4}}{\varpi_{\F}^{2}} \\ 
        y_{3, 1}\varpi_{\F}^{2} & y_{3, 2}\varpi_{\F} & y_{3, 3} & \sfrac{y_{3, 4}}{\varpi_{\F}}\\ 
        y_{4, 1}\varpi_{\F}^{3} & y_{4, 2}\varpi_{\F}^{2} & y_{4, 3}\varpi_{\F} & y_{4, 4}
    \end{pmatrix} = (x_{0} + x_{1}\sigma_{f_{M}} + x_{2}\sigma_{f_{M}}^{2} +x_{3}\sigma_{f_{M}}^{3})z
\end{equation}
where $x_{0} + x_{1}\sigma_{f_{M}} + x_{2}\sigma_{f_{M}}^{2} +x_{3}\sigma_{f_{M}}^{3} \in \mathcal{O}_{\E_{f_{M}}}^{\times}$ is of the form in (\ref{Unramified2}) and $z \in U^{1}_{M}$ is of the form
\[
z = \begin{pmatrix}
       1+ z_{1, 1}\varpi_{\F} & z_{1, 2} & \sfrac{z_{1, 3}}{\varpi_{\F}} & \sfrac{z_{1, 4}}{\varpi_{\F}^{2}} \\ 
        z_{2, 1}\varpi_{\F}^{2} & 1+z_{2, 2}\varpi_{\F} & z_{2, 3} & \sfrac{z_{2, 4}}{\varpi_{\F}} \\ 
        z_{3, 1}\varpi_{\F}^{3} & z_{3, 2}\varpi_{\F}^{2} & 1+z_{3, 3}\varpi_{\F} & z_{3, 4}\\ 
        z_{4, 1}\varpi_{\F}^{4} & z_{4, 2}\varpi_{\F}^{3} & z_{4, 3}\varpi_{\F}^{2} & 1+z_{4, 4}\varpi_{\F}
    \end{pmatrix}
\]
with $y_{i, j}$, $z_{i, j} \in \mathcal{O}_{\F}$.

\begin{lemma}\label{LemmaSupportBiquad}
    Let $g$ be of the form in (\ref{Iwasawa}) and $n^{-}(x)$ be of the form in (\ref{nilpotent}).~Next, set $\alpha_{M} := \sigma_{4}\begin{pmatrix}
    I_{2} & n^{-}(x) \\
    & I_{2}
\end{pmatrix}\begin{pmatrix}
    g & \\
    & g
\end{pmatrix}\sigma_{4}$.~If $\alpha_{M} \in \Su\left(\mathcal{W}_{(f_{M}, \, \chi_{M}, \, \zeta)}\right)$, then 
\[
\alpha_{M} \in \N(4, \F)\varpi_{\F}^{-2}J_{f_{M}} \, \sqcup \, \N(4, \F)\varpi_{\F}^{-1}J_{f_{M}} \, \sqcup \, \N(4, \F)J_{f_{M}}.
\]
\end{lemma}
\begin{proof}
    Let $g$ and $n^{-}(x)$ be as in the hypothesis, then $\alpha_{M}$ is of the form given in (\ref{Shalikaelement}).~Suppose $k >0$.~Then $\alpha_{M} = u\varpi_{\F}^{k}y$ for some $u \in \N(4, \F)$ and $y \in J_{f_{M}}$.~Comparing the $(4, 2)$ and $(4, 4)$-entries of $\alpha_{M}$ and $u\varpi_{\F}^{k}y$, we have that $h_{2, 1}$, $h_{2, 2} \in \mathcal{P}_{\F}$ which implies $h \notin \GL(2, \mathcal{O}_{\F})$.~This is a contradiction and $k \leq 0$.~

Let $k < -2$.~Then $\alpha_{M} = u\varpi_{\F}^{k}y$ for some $u \in \N(4, \F)$ and $y \in J_{f_{M}}$.~Comparing the fourth rows of $\alpha_{M}$ and $u\varpi_{\F}^{k}y$, we have that $y_{4, 1} = y_{4, 3} = 0$ and $y_{4, 2}$, $y_{4, 4} \in \mathcal{P}_{\F}$.~From (\ref{decomposition2}), we see that $y_{4, 1} \in x_{3} +\mathcal{P}_{\F}$ and $y_{4, 3} \in x_{1} +bx_{3} + \mathcal{P}_{\F}$.~Hence, $x_{1}$, $x_{3} \in \mathcal{P}_{\F}$ since $y_{4, 1}$ and $y_{4, 3}$ are zero.~

Using (\ref{decomposition2}) again, we see that $y_{4, 2} \in x_{2} +\mathcal{P}_{\F}$ and $y_{4, 4} \in x_{0} +x_{2}b +\mathcal{P}_{\F}$.~Since both of these entries are contained in $\mathcal{P}_{\F}$, we have that $x_{0}$, $x_{2} \in \mathcal{P}_{\F}$.~This is a contradiction as $x_{0} + x_{1}\sigma_{f_{M}} + x_{2}\sigma_{f_{M}}^{2} +x_{3}\sigma_{f_{M}}^{3} \in \mathcal{O}_{\E_{f_{M}}}^{\times}$ and hence, proves the lemma.
\end{proof}

Before we introduce the following proposition, we note the following.~Let $\text{L}_{f_{M}} = \F\left[\sigma_{f_{M}}^{2}\right]$ denote the unique degree two unramified subextension of $\E_{f_{M}}/\F$ with valuation ring $\mathcal{O}_{\text{L}_{f_{M}}} = \mathcal{O}_{\F}\left[\sigma_{f_{M}}^{2}\right]$.~Then $\psi_{\beta_{f_{M}}}$ is trivial on the subgroup of principal units of $\mathcal{O}_{\text{L}_{f_{M}}}^{\times}$:
\begin{equation}\label{Trivial}
    1 + \mathcal{P}_{\text{L}_{f_{M}}} = \left\{1 + \varpi_{\F}x : x = x_{0} + x_{2}\sigma_{f_{M}}^{2} \in \mathcal{O}_{\text{L}_{f_{M}}}\right\}.
\end{equation}
Furthermore, we may write $x_{0}+x_{2}\sigma_{f_{M}}^{2} \in \mathcal{O}_{\text{L}_{f_{M}}}^{\times}$ as
\begin{equation}\label{ProofDecomposition}
    x_{0}+x_{2}\sigma_{f_{M}}^{2} = (x_{0}+x_{2}\sigma_{f_{M}}^{2})' + \varpi_{\F}\left(w_{0}+w_{2}\sigma_{f_{M}}^{2}\right)
\end{equation}
where $(x_{0}+x_{2}\sigma_{f_{M}}^{2})' \in \mu_{\text{L}_{f_{M}}}'$ and $w_{0}+w_{2}\sigma_{f_{M}}^{2} \in \mathcal{O}_{\text{L}_{f_{M}}}$.~
\begin{proposition}\label{twistedbiquad}
    Let $\pi_{(f_{M}, \, \chi_{M}, \, \zeta)}$ be a biquadratic supercuspidal representation of $\GL(4, \F)$. Then $L\left(s, \pi_{(f_{M}, \, \chi_{M}, \, \zeta)}, \wedge^{2}\right)$ has a pole at $s = 0$ if and only if $\omega_{(f_{M}, \, \chi_{M}, \, \zeta)}$ is trivial.
\end{proposition}
\begin{proof}
        If $L\left(s, \pi_{(f_{M}, \, \chi_{M}, \, \zeta)}, \wedge^{2}\right)$ has a pole at $s = 0$, then $\omega_{(f_{M}, \, \chi_{M}, \, \zeta)}$ is necessarily trivial from Theorem \ref{Jo}.~To prove the converse, suppose $\omega_{(f_{M}, \, \chi_{M}, \, \zeta)}$ is trivial.~We claim that $\Lambda_{0}$ evaluated at $\pi_{(f_{M}, \, \chi_{M}, \, \zeta)}(\sigma_{4})\mathcal{W}_{(f_{M}, \, \chi_{M}, \, \zeta)}$ is non-zero.
    
    Using our coset representative choices in (\ref{nilpotent}) and (\ref{Iwasawa}), we have that 
    \begin{align*}
        &\Lambda_{0}\left(\pi_{(f_{M}, \, \chi_{M}, \, \zeta)}(\sigma_{4})\mathcal{W}_{(f_{M}, \, \chi_{M}, \, \zeta)}\right) 
         \\ & = \int\displaylimits_{\F^{\times}\N(2, \, \F) \backslash \GL(2, \, \F)}\int\displaylimits_{\F}\mathcal{W}_{(f_{M}, \, \chi_{M}, \, \zeta)}\left(\sigma_{4}\begin{pmatrix}
    I_{2} & n^{-}(x) \\
    & I_{2}
\end{pmatrix}\begin{pmatrix}
    g & \\
    & g
\end{pmatrix}\sigma_{4}\right)\, dx \, d^{\times}g.
    \end{align*}
Furthermore, Lemma \ref{LemmaSupportBiquad} tells us that we need only consider the double cosets 
\[
\F^{\times}\N(2, \F)\begin{pmatrix}
    \varpi_{\F}^{-2k} & \\
    & 1
\end{pmatrix}\GL(2, \mathcal{O}_{\F})
\]
with $0 \leq k \leq 2$ for the overall value of $\Lambda_{0}\left(\pi_{(f_{M}, \, \chi_{M}, \, \zeta)}(\sigma_{4})\mathcal{W}_{(f_{M}, \, \chi_{M}, \, \zeta)}\right)$.~Furthermore, for each $0 \leq k \leq 2$, if an element $\alpha$ of the form in (\ref{Shalikaelement}) with respect to the double coset $\F^{\times}\N(2, \F)\begin{pmatrix}
    \varpi_{\F}^{-2k} & \\
    & 1
\end{pmatrix}\GL(2, \mathcal{O}_{\F})$ is contained inside $\Su\left(\mathcal{W}_{(f_{M}, \, \chi_{M}, \, \zeta)}\right)$, then $\alpha$ is contained inside $\N(4, \F)\varpi_{\F}^{-k}J_{f_{M}}$ as $\det(\alpha) \in \mathcal{P}_{\F}^{-4k}$.~Hence, we may consider three separate integrals $T_{0}$, $T_{1}$, and $T_{2}$ such that 
\begin{equation}\label{ThreeIntegrals}
    \Lambda_{0}\left(\pi_{(f_{M}, \, \chi_{M}, \, \zeta)}(\sigma_{4})\mathcal{W}_{(f_{M}, \, \chi_{M}, \, \zeta)}\right) =T_{0} + T_{1}+T_{2}
\end{equation}
where 
\begin{align*}
   & T_{k} \\
   & = \int\displaylimits_{\F^{\times}\N(2, \, \F) \backslash \F^{\times}\N(2, \, \F)\begin{pmatrix}
    \varpi_{\F}^{-2k} & \\
    & 1
\end{pmatrix}\GL(2, \, \mathcal{O}_{\F})}\int\displaylimits_{\F}\mathcal{W}_{(f_{M}, \, \chi_{M}, \, \zeta)}\left(\sigma_{4}\begin{pmatrix}
    I_{2} & n^{-}(x) \\
    & I_{2}
\end{pmatrix}\begin{pmatrix}
    g & \\
    & g
\end{pmatrix}\sigma_{4}\right)\, dx \, d^{\times}g
\end{align*}
for $0 \leq k \leq 2$.

We first compute $T_{0}$.~If an element $\alpha$ of the form in (\ref{Shalikaelement}) is contained inside $\N(4, \F)J_{f_{M}}$, then $\alpha = uy$ for some $u = (u_{i, j}) \in \N(4, \F)$ and $y=(x_{0} + x_{1}\sigma_{f_{M}} + x_{2}\sigma_{f_{M}}^{2} +x_{3}\sigma_{f_{M}}^{3})z \in J_{f_{M}}$.~Since $\det(\alpha) \in \mathcal{O}_{\F}^{\times}$, we may choose representatives of the form in (\ref{Iwasawa}) such that $r = 0$.~Comparing the $(4, 1)$ and $(4, 3)$-entries of $\alpha$ and $uy$, we have that $y_{4, 1} = y_{4, 3} = 0$ which implies $x_{1}$, $x_{3} \in \mathcal{P}_{\F}$ as $y_{4, 1} \in x_{3} +\mathcal{P}_{\F}$ and $y_{4, 3} \in x_{1}+bx_{3} +\mathcal{P}_{\F}$ from (\ref{decomposition2}).

Next, we compare the $(3, 1)$ and $(4, 2)$-entries of $\alpha$ and $uy$ to see that $h_{2, 1} \in (x_{2} + \mathcal{P}_{\F})\varpi_{\F}^{2}$ from (\ref{decomposition2}); this implies $h \in I$ and from (\ref{IwahoriReps}), we may choose representatives of the form 
\[
\begin{pmatrix}
    h_{1, 1} & \\
    h_{2, 1} & h_{2, 2}
\end{pmatrix}
\]
where $h_{1, 1}$, $h_{2, 2} \in \mathcal{O}_{\F}^{\times}$.~Comparing the $(3, 3)$ and $(4, 4)$-entries of $\alpha$ and $uy$, we have that $h_{2, 2} \in x_{0}+x_{2}b + \mathcal{P}_{\F}$ from (\ref{decomposition2}).~

To see where $x$ lies, we first suppose $h_{2, 1}$ is zero.~Comparing the $(3, 2)$-entries of $\alpha$ and $uy$, we have that $x \in \mathcal{P}_{\F}^{2}$ as  $y_{3, 2} \in \mathcal{P}_{\F}$ from $(\ref{decomposition2})$ and $h_{1, 1} \in \mathcal{O}_{\F}^{\times}$.~Now suppose $h_{2, 1}$ is non-zero.~Comparing the $(3, 2)$ and $(3, 4)$-entries of $\alpha$ and $uy$, we get the equation
\begin{equation}\label{TB1}
    u_{3, 4} = \frac{h_{1, 1}x-y_{3, 2}\varpi_{\F}}{h_{2, 1}} = -\frac{y_{3, 4}}{h_{2, 2}\varpi_{\F}}.
\end{equation}
Hence, equation (\ref{TB1}) implies 
\[
x = \frac{(h_{2, 2}y_{3, 2}-h_{2, 1}y_{3, 4}\varpi_{\F}^{-2})\varpi_{\F}}{h_{1, 1}h_{2, 2}} \in \mathcal{P}_{\F}^{2}
\]
as $h_{2, 2}y_{3, 2}-h_{2, 1}y_{3, 4}\varpi_{\F}^{-2} \in \mathcal{P}_{\F}$ from (\ref{decomposition2}).

To see where $h_{1, 1}$ lies, we first compare the $(2, 3)$-entries of $\alpha$ and $uy$ and see that $u_{2, 3} \in \mathcal{O}_{\F}$ as $y_{2, 3} \in \mathcal{P}_{\F}$ from (\ref{decomposition2}).~Next, we compare the $(2, 2)$ and $(2, 4)$-entries of $\alpha$ and $uy$ to get the equation
\begin{equation}\label{TB2}
    u_{2, 4} = \frac{h_{1, 1}-y_{2, 2}-u_{2, 3}y_{3, 2}\varpi_{\F}}{h_{2, 1}} = -\frac{y_{2, 4}+u_{2, 3}y_{3, 4}}{h_{2, 2}\varpi_{\F}^{2}}.
\end{equation}
Using equations (\ref{TB2}) and (\ref{decomposition2}), we see that
\[
h_{1, 1} \in \frac{x_{0}^{2}+x_{0}x_{2}b-x_{2}^{2}a}{x_{0}+x_{2}b} + \mathcal{P}_{\F}.
\]
Hence, $\alpha$ lies inside
\[
\begin{pmatrix}
    \frac{x_{0}^{2}+x_{0}x_{2}b-x_{2}^{2}a}{x_{0}+x_{2}b} + \mathcal{P}_{\F} & & & \\
    & \frac{x_{0}^{2}+x_{0}x_{2}b-x_{2}^{2}a}{x_{0}+x_{2}b} + \mathcal{P}_{\F} & & \\
    (x_{2} + \mathcal{P}_{\F})\varpi_{\F}^{2} & \mathcal{P}_{\F}^{2} & x_{0}+x_{2}b + \mathcal{P}_{\F} & \\
    & (x_{2} + \mathcal{P}_{\F})\varpi_{\F}^{2} & & x_{0}+x_{2}b + \mathcal{P}_{\F}
\end{pmatrix}.
\]

Our computations above imply that we may decompose $\alpha$ as follows:
\[
\alpha = u (x_{0}+x_{2}\sigma_{f_{M}}^{2})z,
\]
where the entries of $u$ satisfy
\[
u_{1, 3} = u_{2, 4} = \frac{h_{1, 1}a}{h_{2, 1}a+h_{2, 2}b\varpi_{\F}^{2}}, \; \; u_{1, 2} = u_{1, 4} = u_{2, 3} = u_{3, 4} = 0
\]
and the entries of $z$ satisfy
\begin{align*}
     & 1+ z_{1, 1}\varpi_{\F} = 1+ z_{2, 2}\varpi_{\F} = \frac{h_{1, 1}h_{2, 2}x_{0}(x_{0}+x_{2}b)\varpi_{\F}^{4} + h_{2, 1}x_{2}a(h_{2, 1}x_{2}a-h_{2, 2}x_{0}\varpi_{\F}^{2})}{(x_{0}^{2}+x_{0}x_{2}b-x_{2}^{2}a)(h_{2, 2}x_{0}\varpi_{\F}^{2}-h_{2, 1}x_{2}a)\varpi_{\F}^{2}} \\
     & z_{1, 2} = \frac{h_{1, 1}xx_{2}a(h_{2, 1}x_{2}a+h_{1, 1}(x_{0}+x_{2}b)\varpi_{\F}^{2}-h_{2, 2}x_{0}\varpi_{\F}^{2})}{(x_{0}^{2}+x_{0}x_{2}b-x_{2}^{2}a)(h_{2, 2}x_{0}\varpi_{\F}^{2}-h_{2, 1}x_{2}a)\varpi_{\F}^{2}} \\
     & z_{1, 3}\varpi_{\F}^{-1} = z_{2, 4}\varpi_{\F}^{-1} = \frac{h_{2, 2}x_{2}a(h_{2, 1}x_{2}a+h_{1, 1}(x_{0}+x_{2}b)\varpi_{\F}^{2}-h_{2, 2}x_{0}\varpi_{\F}^{2})}{(x_{0}^{2}+x_{0}x_{2}b-x_{2}^{2}a)(h_{2, 2}x_{0}\varpi_{\F}^{2}-h_{2, 1}x_{2}a)\varpi_{\F}^{2}} \\
     & z_{3, 1}\varpi_{\F}^{3} = z_{4, 2}\varpi_{\F}^{3} = \frac{x_{0}(h_{2, 1}h_{2, 2}x_{0}\varpi_{\F}^{2}-h_{2, 1}^{2}ax_{2}-h_{1, 1}h_{2, 2}x_{2}\varpi_{\F}^{4})}{(x_{0}^{2}+x_{0}x_{2}b-x_{2}^{2}a)(h_{2, 2}x_{0}\varpi_{\F}^{2}-h_{2, 1}x_{2}a)} \\
     & z_{3, 2}\varpi_{\F}^{2} = \frac{h_{1, 1}x(h_{2, 2}x_{0}^{2}\varpi_{\F}^{2}-x_{2}a(h_{2, 1}x_{0}+h_{1, 1}x_{2}\varpi_{\F}^{2}))}{(x_{0}^{2}+x_{0}x_{2}b-x_{2}^{2}a)(h_{2, 2}x_{0}\varpi_{\F}^{2}-h_{2, 1}x_{2}a)} \\
     & 1+ z_{3, 3}\varpi_{\F} = 1+ z_{4, 4}\varpi_{\F} = \frac{h_{2, 2}(h_{2, 2}x_{0}^{2}\varpi_{\F}^{2}-x_{2}a(h_{2, 1}x_{0}+h_{1, 1}x_{2}\varpi_{\F}^{2}))}{(x_{0}^{2}+x_{0}x_{2}b-x_{2}^{2}a)(h_{2, 2}x_{0}\varpi_{\F}^{2}-h_{2, 1}x_{2}a)}
\end{align*}
with $z_{i, j} = 0$ for all other $i$ and $j$.~For the entries of $z$, we note that $x_{0}^{2}+x_{0}x_{2}b-x_{2}^{2}a \in \mathcal{O}_{\F}^{\times}$ as $\det(x_{0}+x_{2}\sigma_{f_{M}}^{2}) = (x_{0}^{2}+x_{0}x_{2}b-x_{2}^{2}a)^{2} \in \mathcal{O}_{\F}^{\times}$.~

Since $\psi_{4}(u) = \psi_{\beta_{f_{M}}}(z) = 1$, we have by definition of $\mathcal{W}_{(f_{M}, \, \chi_{M}, \, \zeta)}$ and equations (\ref{Trivial}) and (\ref{ProofDecomposition}) that 
\begin{equation}\label{TB3}
    \mathcal{W}_{(f_{M}, \, \chi_{M}, \, \zeta)}(\alpha) =\varphi_{M}\left(x_{0}+x_{2}\sigma_{f_{M}}^{2}\right) = \chi_{M}\left((x_{0}+x_{2}\sigma_{f_{M}}^{2})'\right).
\end{equation}
Hence, equation (\ref{TB3}) implies that $T_{0}$ is zero if $\chi_{M}$ is non-trivial on $\mathcal{O}_{\text{L}_{f_{M}}}^{\times}$ as we are integrating over a compact subgroup $I$ and $\mathcal{P}_{\F}^{2}$.~If $\chi_{M}$ is trivial on $\mathcal{O}_{\text{L}_{f_{M}}}^{\times}$, then $T_{0}$ is equal to a positive real volume.

Next, we compute $T_{1}$.~If $\alpha$ of the form in (\ref{Shalikaelement}) is contained inside $\N(4, \F)\varpi_{\F}^{-1}J_{f_{M}}$, then $\alpha = u\varpi_{\F}^{-1}y$ for some $u = (u_{i, j}) \in \N(4, \F)$ and $y=(x_{0} + x_{1}\sigma_{f_{M}} + x_{2}\sigma_{f_{M}}^{2} +x_{3}\sigma_{f_{M}}^{3})z \in J_{f_{M}}$.~Since $\det(\alpha) \in \mathcal{P}_{\F}^{-4}$, we may choose representatives of the form in (\ref{Iwasawa}) such that $r = -2$.~Comparing the $(4, 1)$ and $(4, 3)$-entries of $\alpha$ and $u\varpi_{\F}^{-1}y$, we have that $y_{4, 1} = y_{4, 3} = 0$ which implies $x_{1}$, $x_{3} \in \mathcal{P}_{\F}$ as $y_{4, 1} \in x_{3} +\mathcal{P}_{\F}$ and $y_{4, 3} \in x_{1}+bx_{3} +\mathcal{P}_{\F}$ from (\ref{decomposition2}).

We compare the $(3, 1)$ and $(4, 2)$-entries of $\alpha$ and $u\varpi_{\F}^{-1}y$ to see that $h_{2, 1} \in (x_{2} + \mathcal{P}_{\F})\varpi_{\F}$; this implies $h \in I$ and from (\ref{IwahoriReps}), we may choose representatives of the form 
\[
\begin{pmatrix}
    \varpi_{\F}^{-2} & \\
    & 1
\end{pmatrix}\begin{pmatrix}
    h_{1, 1} & \\
    h_{2, 1} & h_{2, 2}
\end{pmatrix}
\]
where $h_{1, 1}$, $h_{2, 2} \in \mathcal{O}_{\F}^{\times}$.~Comparing the $(3, 3)$ and $(4, 4)$-entries of $\alpha$ and $u\varpi_{\F}^{-1}y$, we see that $h_{2, 2} = y_{3, 3}\varpi_{\F}^{-1}$.~This implies $x_{0} \equiv -x_{2}b \imod{\mathcal{P}_{\F}}$ from (\ref{decomposition2}) and $x_{2} \in \mathcal{O}_{\F}^{\times}$.~Let $h_{2, 2} \in v + \mathcal{P}_{\F}$ for some $v \in \mu_{\F}'$.~

To see where $x$ lies, we compare the $(3, 2)$ and $(3, 4)$-entries of $\alpha$ and $u\varpi_{\F}^{-1}y$.~In doing so, we get the equation
\begin{equation}\label{TB4}
    u_{3, 4} = \frac{h_{1, 1}x-y_{3, 2}\varpi_{\F}^{2}}{h_{2, 1}\varpi_{\F}^{2}} = -\frac{y_{3, 4}}{h_{2, 2}\varpi_{\F}^{2}}.
\end{equation}
Equation (\ref{TB4}) implies 
\[
x = \frac{h_{2, 2}y_{3, 2}\varpi_{\F}^{2}-h_{2, 1}y_{3, 4}}{h_{1, 1}h_{2, 2}} \in \mathcal{P}_{\F}^{2}
\]
as $y_{3, 4} \in \mathcal{P}_{\F}$ from (\ref{decomposition2}).

To see where $h_{1, 1}$ lies, we first compare the $(2, 1)$-entries of $\alpha$ and $u\varpi_{\F}^{-1}y$ and see that $u_{2, 3} \in \mathcal{O}_{\F}$ as $y_{2, 1} \in \mathcal{P}_{\F}$ from (\ref{decomposition2}).~Next, we compare the $(2, 2)$ and $(2, 4)$-entries of $\alpha$ and $u\varpi_{\F}^{-1}y$ to get the equation
\begin{equation}\label{TB5}
   u_{2, 4} = \frac{h_{1, 1}-y_{2, 2}\varpi_{\F}-u_{2, 3}y_{3, 2}\varpi_{\F}^{2}}{h_{2, 1}\varpi_{\F}^{2}} = -\frac{y_{2, 4}+u_{2, 3}y_{3, 4}\varpi_{\F}}{h_{2, 2}\varpi_{\F}^{3}}.~
\end{equation}
Using equations (\ref{TB5}) and (\ref{decomposition2}), we see that 
\[
h_{1, 1} \in -\frac{x_{2}^{2}a}{v} + \mathcal{P}_{\F}.
\]
Hence, $\alpha$ lies inside
\[
\begin{pmatrix}
    \left(-\frac{x_{2}^{2}a}{v} + \mathcal{P}_{\F}\right)\varpi_{\F}^{-2} & & & \\
    & \left(-\frac{x_{2}^{2}a}{v} + \mathcal{P}_{\F}\right)\varpi_{\F}^{-2} & & \\
    (x_{2} + \mathcal{P}_{\F})\varpi_{\F} & \mathcal{O}_{\F} & v + \mathcal{P}_{\F} & \\
    & (x_{2} + \mathcal{P}_{\F})\varpi_{\F} & & v + \mathcal{P}_{\F}
\end{pmatrix}.
\]

Our computations above imply that we may decompose $\alpha$ as follows:
\[
\alpha = u \varpi_{\F}^{-1}x_{2}\left(-b+\sigma_{f_{M}}^{2}\right)z,
\]
where the entries of $u$ satisfy
\begin{align*}
     & u_{1, 2} = -\frac{x(h_{1, 1}+x_{2}b\varpi_{\F})}{h_{2, 1}\varpi_{\F}^{2}}, \; \; 
     u_{1, 3} = u_{2, 4} = \frac{h_{1, 1}+x_{2}b\varpi_{\F}}{h_{2, 1}\varpi_{\F}^{2}} \\
     & u_{1, 4} = -\frac{xb(h_{1, 1}+x_{2}b\varpi_{\F})}{h_{2, 1}\varpi_{\F}^{4}}, \; \; 
     u_{2, 3} = 0, \; \; 
     u_{3, 4} = \frac{x(h_{2, 1}b-(h_{1, 1}+x_{2}b\varpi_{\F}))}{h_{2, 1}\varpi_{\F}^{2}}
\end{align*}
and the entries of $z$ satisfy
\begin{align*}
    & 1 + z_{1, 1}\varpi_{\F} = 1 + z_{2, 2}\varpi_{\F} = \frac{h_{2, 1}}{x_{2}\varpi_{\F}}, \; \; 
    z_{1, 2} = \frac{xb(h_{2, 1}-x_{2}\varpi_{\F})}{x_{2}\varpi_{\F}^{3}}, \; \;  
    z_{1, 3}\varpi_{\F}^{-1} = z_{2, 4}\varpi_{\F}^{-1} = \frac{h_{2, 2}}{x_{2}\varpi_{\F}}  \\
    & z_{1, 4}\varpi_{\F}^{-2}= \frac{h_{2, 2}x(h_{2, 1}b-(h_{1, 1}+x_{2}b\varpi_{\F}))}{h_{2, 1}x_{2}\varpi_{\F}^{3}}, \; \; 
    z_{3, 1}\varpi_{\F}^{3} = z_{4, 2}\varpi_{\F}^{3} = \frac{b(h_{2, 1}-x_{2}\varpi_{\F})\varpi_{\F}}{x_{2}a} \\
    & z_{3, 2}\varpi_{\F}^{2} = \frac{xb^{2}(h_{2, 1}-x_{2}\varpi_{\F})}{x_{2}a\varpi_{\F}}, \; \; z_{3, 4} = \frac{h_{2, 2}xb(h_{2, 1}b-(h_{1, 1}+x_{2}b\varpi_{\F}))}{h_{2, 1}x_{2}a\varpi_{\F}} \\
    & 1 + z_{3, 3}\varpi_{\F} = 1 + z_{4, 4}\varpi_{\F} = \frac{h_{2, 2}(h_{2, 1}b-(h_{1, 1}+x_{2}b\varpi_{\F}))\varpi_{\F}}{h_{2, 1}x_{2}a}
\end{align*}
with $z_{i, j} = 0$ for all other $i$ and $j$.~

Since $\omega_{(f_{M}, \, \chi_{M}, \, \zeta)}$ is trivial, $\psi_{4}(u) = \psi_{\F}(bx\varpi_{\F}^{-2})^{-1}$, and $\psi_{\beta_{f_{M}}}(z) = \psi_{\F}(bx\varpi_{\F}^{-2})$, we have by definition of $\mathcal{W}_{(f_{M}, \, \chi_{M}, \, \zeta)}$ and equations (\ref{Trivial}) and (\ref{ProofDecomposition}) that
\begin{equation}\label{TB6}
    \mathcal{W}_{(f_{M}, \, \chi_{M}, \, \zeta)}(\alpha) =\varphi_{M}\left(-b+\sigma_{f_{M}}^{2}\right) = \chi_{M}\left(\left(-b+\sigma_{f_{M}}^{2}\right)'\right).~
\end{equation}
Hence, equation (\ref{TB6}) tells us that $T_{1}$ is equal to $\chi_{M}\left(\left(-b+\sigma_{f_{M}}^{2}\right)'\right) \text{V}_{1}$ for some positive real volume $V_{1}$ as we are integrating over a compact subgroup $I$ and $\mathcal{P}_{\F}^{2}$.

Lastly, we compute $T_{2}$.~If $\alpha$ of the form in (\ref{Shalikaelement}) is contained inside $\N(4, \F)\varpi_{\F}^{-2}J_{f_{M}}$, then $\alpha = u\varpi_{\F}^{-2}y$ for some $u = (u_{i, j}) \in \N(4, \F)$ and $y=(x_{0} + x_{1}\sigma_{f_{M}} + x_{2}\sigma_{f_{M}}^{2} +x_{3}\sigma_{f_{M}}^{3})z \in J_{f_{M}}$.~Since $\det(\alpha) \in \mathcal{P}_{\F}^{-8}$, we may choose representatives of the form in (\ref{Iwasawa}) such that $r = -4$.~Comparing the $(4, 1)$ and $(4, 3)$-entries of $\alpha$ and $u\varpi_{\F}^{-2}y$, we have that $y_{4, 1} = y_{4, 3} = 0$ which implies $x_{1}$, $x_{3} \in \mathcal{P}_{\F}$ as $y_{4, 1} \in x_{3} + \mathcal{P}_{\F}$ and $y_{4, 3} \in x_{1}+bx_{3} +\mathcal{P}_{\F}$ from (\ref{decomposition2}).

Next, we compare the $(3, 3)$ and $(4, 4)$-entries of $\alpha$ and $u\varpi_{\F}^{-2}y$ and see that $h_{2, 2} = y_{3, 3}\varpi_{\F}^{-2}$.~This implies $x_{0} \equiv -x_{2}b \imod{\mathcal{P}_{\F}}$ from (\ref{decomposition2}) and $x_{2} \in \mathcal{O}_{\F}^{\times}$.~Comparing the $(3, 1)$ and $(4, 2)$-entries of $\alpha$ and $u\varpi_{\F}^{-2}y$, we have that $h_{2, 1} \in x_{2}+\mathcal{P}_{\F}$; this implies $h \in IsI$ and from (\ref{IwahoriReps2}), we may choose representatives of the form 
\[
\begin{pmatrix}
    \varpi_{\F}^{-4} & \\
    & 1
\end{pmatrix}\begin{pmatrix}
     & h_{1, 2}\\
    h_{2, 1} & h_{2, 2}
\end{pmatrix}
\]
where $h_{1, 2}$, $h_{2, 1} \in \mathcal{O}_{\F}^{\times}$.~

To see where $x$ lies, we compare the $(3, 2)$ and $(3, 4)$-entries of $\alpha$ and $u\varpi_{\F}^{-2}y$, we get the equation
\begin{equation}\label{TB7}
    u_{3, 4} = \frac{h_{1, 2}x-y_{3, 4}\varpi_{\F}}{h_{2, 2}\varpi_{\F}^{4}} = -\frac{y_{3, 2}}{h_{2, 1}\varpi_{\F}}.
\end{equation}
Equation (\ref{TB7}) implies 
\[
x = -\frac{(h_{2, 2}y_{3, 2}\varpi_{\F}^{2}-h_{2, 1}y_{3, 4})\varpi_{\F}}{h_{1, 2}h_{2, 1}} \in \mathcal{P}_{\F}^{2}
\]
as $y_{3, 4} \in \mathcal{P}_{\F}$ from (\ref{decomposition2}).

To see where $h_{1, 2}$ lies, we first compare the $(2, 1)$-entries of $\alpha$ and $u\varpi_{\F}^{-2}y$ and see that $u_{2, 3} \in \mathcal{O}_{\F}$ as $y_{2, 1} \in \mathcal{P}_{\F}$ from (\ref{decomposition2}).~Next, we compare the $(2, 2)$ and $(2, 4)$-entries of $\alpha$ and $u\varpi_{\F}^{-2}y$ to get the equation
\begin{equation}\label{TB8}
    u_{2, 4} = \frac{h_{1, 2}-y_{2, 4}-u_{2, 3}y_{3, 4}\varpi_{\F}}{h_{2, 2}\varpi_{\F}^{4}} = -\frac{y_{2, 2}+u_{2, 3}y_{3, 2}\varpi_{\F}}{h_{2, 1}\varpi_{\F}^{2}}.
\end{equation}
Using equations (\ref{TB8}) and (\ref{decomposition2}), we see that $h_{1, 2} \in x_{2}a + \mathcal{P}_{\F}$.~Hence, $\alpha$ lies inside
\[
\begin{pmatrix}
     & & (x_{2}a + \mathcal{P}_{\F})\varpi_{\F}^{-4} & \\
    &  & & (x_{2}a + \mathcal{P}_{\F})\varpi_{\F}^{-4} \\
    x_{2} + \mathcal{P}_{\F} &  & \mathcal{O}_{\F} & \mathcal{P}_{\F}^{-2} \\
    & x_{2} + \mathcal{P}_{\F} & & \mathcal{O}_{\F}
\end{pmatrix}.
\]

Our computations above imply that we may decompose $\alpha$ as follows:
\[
\alpha = u \varpi_{\F}^{-2}x_{2}\left(-b+\sigma_{f_{M}}^{2}\right)z,
\]
where the entries of $u$ satisfy
\[
    u_{1, 3} = u_{2, 4} = \frac{x_{2}b}{h_{2, 1}\varpi_{\F}^{2}}, \; \; u_{1, 2} = u_{1, 4} = u_{2, 3} = u_{3, 4} = 0
\]
and the entries of $z$ satisfy
\begin{align*}
    & 1+z_{1, 1}\varpi_{\F} = 1+z_{2, 2}\varpi_{\F} = \frac{h_{2, 1}}{x_{2}}, \; \; z_{1, 3}\varpi_{\F}^{-1} = z_{2, 4}\varpi_{\F}^{-1} = \frac{h_{2, 2}}{x_{2}} \\
    & z_{1, 4}\varpi_{\F}^{-2} = \frac{h_{1, 2}x}{x_{2}\varpi_{\F}^{4}}, \; \; z_{3, 1}\varpi_{\F}^{3} = z_{4, 2}\varpi_{\F}^{3} = \frac{b(h_{2, 1}-x_{2})\varpi_{\F}^{2}}{x_{2}a}, \; \; z_{3, 4} = \frac{h_{1, 2}xb(h_{2, 1}-x_{2})}{x_{2}ah_{2, 1}\varpi_{\F}^{2}} \\
    & 1+z_{3, 3}\varpi_{\F} = 1+z_{4, 4}\varpi_{\F} = \frac{h_{1, 2}h_{2, 1}+h_{2, 2}b(h_{2, 1}-x_{2})\varpi_{\F}^{2}}{h_{2, 1}x_{2}a}
\end{align*}
with $z_{i, j} = 0$ for all other $i$ and $j$.

Since $\omega_{(f_{M}, \, \chi_{M}, \, \zeta)}$ is trivial and $\psi_{4}(u) = \psi_{\beta_{f_{M}}}(z) = 1$, we have by definition of $\mathcal{W}_{(f_{M}, \, \chi_{M}, \, \zeta)}$ and equations (\ref{Trivial}) and (\ref{ProofDecomposition}) that
\begin{equation}\label{TB9}
    \mathcal{W}_{(f_{M}, \, \chi_{M}, \, \zeta)}(\alpha) =\varphi_{M}\left(-b+\sigma_{f_{M}}^{2}\right) = \chi_{M}\left(\left(-b+\sigma_{f_{M}}^{2}\right)'\right).~
\end{equation}
Hence, equation (\ref{TB9}) tells us that $T_{2}$ is equal to $\chi_{M}\left(\left(-b+\sigma_{f_{M}}^{2}\right)'\right)\text{V}_{2}$ for some positive real volume $V_{2}$ as we are integrating over a compact subgroup $IsI$ and $\mathcal{P}_{\F}^{2}$.

Our calculations show that in either situation with $\chi_{M}$ trivial or non-trivial, $\Lambda_{0}$ evaluated at $\pi_{(f_{M}, \, \chi_{M}, \, \zeta)}(\sigma_{4})\mathcal{W}_{(f_{M}, \, \chi_{M}, \, \zeta)}$ is non-zero from (\ref{ThreeIntegrals}).~The proof of the converse now follows from Theorem \ref{Jo}.
\qedhere
\end{proof}

\subsection{Proof of Theorem \ref{Main}}\label{End}

We are now ready to prove Theorem \ref{Main}.
\begin{proof}[Proof of Theorem \ref{Main}]
By Theorem \ref{Jiang}, statements (1) and (2) are equivalent, and each is equivalent to $L(s,\pi,\wedge^{2})$ having a pole at $s=0$.~It therefore remains to determine, in each of the four cases under consideration, when this pole occurs.~

If $\pi$ is a depth-zero supercuspidal representation, then \cite[Theorems 3.14 and 3.15]{YZ1} shows that $L(s,\pi,\wedge^{2})$ has a pole at $s=0$ if and only if $\pi$ has trivial central character.~If $\pi=\pi_{(\overline{v},\,\varphi,\,\zeta)}$ is a simple supercuspidal representation, then \cite[Theorem 3.5 (2)]{YZ2} shows that $L(s,\pi,\wedge^{2})$ has a pole at $s=0$ if and only if $\pi$ has trivial central character and $\zeta=\pm1$.~

Finally, if $\pi$ is a middle or biquadratic supercuspidal representation, then Propositions \ref{middle} and \ref{twistedbiquad}, respectively, show that $L(s,\pi,\wedge^{2})$ has a pole at $s=0$ if and only if $\pi$ has trivial central character.~Combining these results with Theorem \ref{Jiang} proves all the asserted equivalences.~
\end{proof}

\begin{remark}\label{RemarkMain}
    \normalfont
    Assume now that $\F$ has odd residual characteristic $p$.~We note that a depth-zero, simple, middle, or biquadratic supercuspidal representation has a tamely ramified central character $\omega$, i.e.~$\omega$ is trivial on $1 + \mathcal{P}_{\F}$.
    Let $\Psi = \eta \, \circ \, \text{det}$ be a quasi-character of $\GL(4, \F)$ where $\eta$ is a quasi-character of $\F^{\times}$.
    
    Suppose $\pi \otimes \Psi$ has trivial central character, then $\eta^{4}$ is tamely ramified as the central character of $\pi$ is tamely ramified.~Since $p$ is odd, we have that $\eta$ is tamely ramified.~This implies that $\pi \otimes \Psi$ is of the same type as $\pi$, and hence is covered by the cases in Theorem \ref{Main}.~
\end{remark}

\newcommand{\etalchar}[1]{$^{#1}$}

\end{document}